\documentclass[12pt]{article} 
\usepackage{amsmath, amssymb}
\usepackage{amsthm} 
\usepackage[all]{xy}
\theoremstyle{plain}

\theoremstyle{definition}

\def\I{{\cal I}}
\def\bgn{\begin}

\def\CL{\text{\rm CL}}

\def\J{{\cal J}}

\def\Ob{\text{\rm Ob}}
\def\1{{[1]}}
\def\2{{[2]}}
\def\3{{[3]}}
\def\({\left(}
\def\){\right)}
\def\s-circ{\,{\scriptstyle{\circ}}\,}
\def\<<{<\negthinspace \negthinspace<}
\def\Ad{\text{\rm Ad}}
\def\ad{\text{\rm ad}}

\def\bgn{\begin}

\def\endaln{\end{align}}

\def\<{<\negthinspace \negthinspace <}

\def\({\left(}
\def\){\right)}

\def\[{\big[\neg\big[}
\def\]{\big]\neg\big]}
\def\al{\al}

\def\M{{\cal M}}

\def\a{\alpha}
\def\b{\beta}

\def\e{\varepsilon}
\def\gam{\gamma}
\def\Gam{\Gamma}

\def\del{\delta}
\def\lam{\lambda}
\def\ome{\omega}
\def\Ome{\Omega}
\def\sig{\sigma}
\def\Sig{\Sigma}

\def\Q{\Bbb Q}
\def\R{\Bbb R}
\def\C{\Bbb C}
\def\H{\text{\rm H}}

\def\O{\Bbb O}

\def\M{\frak M}

\def\w{\wedge}
\def\({\left(}
\def\){\right)}

\def\neg{\negthinspace}
\def\h{\hat}

\def\til{\tilde}
\def\wtil{\widetilde}
\def\ch{\check}

\def\ol{\overline}
\def\pa{\partial}

\def\ran{\rangle} 
\def\lan{\langle}

\def\trian{\triangle}

\def\bsh{\backslash}

\def\:{\, :\,}
\def\CL{\text{\rm CL}}
\def\TT{T\oplus T^*}
\def\complex{generalized complex }
\def\K\"ahler{generalized K\"ahler }

\def\10{\displaystyle L^{10}}
\def\2{\displaystyle L^2}
\def\c0{\displaystyle C^0}

\def\10{\displaystyle L^{10}}
\def\2{\displaystyle L^2}
\def\del{\delta}
\def\del2{\displaystyle L^2_{0,\delta}}
\def\c0{\displaystyle C^0}

\def\del{\delta}

\def\K{{\cal K}}
\def\M-A{\text{\rm Monge-Amp\`ere}}
\def\O{{\cal O}}

\def\M-A{\text{\rm Monge-Amp\`ere}}
\def\[{\big[\,}
\def\]{\,\big]}
\def\MC{\text{\rm Maurer-Cartan equation}}
\def\K\"ahler{generalized K\"ahler}
\def\oL{\ol L}
\makeatletter
\def\address#1#2{\begingroup
\noindent\parbox[t]{7.8cm}{%
\small{\scshape\ignorespaces#1}\par\vskip1ex
\noindent\small{\itshape E-mail address}%
\/: #2\par\vskip4ex}\hfill%
\endgroup}%
\makeatother
\title{{Deformations of Generalized K\"ahler Structures and Bihermitian Structures}} 
\author{  
\textsc{Ryushi Goto$^{*}$} 
}
\date{2009/10/9} 
\begin{document}
\maketitle
\footnote{ 
\indent\,\,\, deformation theory, \MC 
}
\footnote{ 
$^{*}$Partly supported by the Grant-in-Aid for Scientific Research (C),
Japan Society for the Promotion \par \indent \,\,\,\,\,of Science. 

}
\begin{abstract}
\noindent
Let $(X,J)$ be a compact K\"ahler manifold with a non-zero holomorphic Poisson structure $\b$.
If the obstruction space for deformations of \complex structures 
on $(X, J)$ vanishes, 
we obtain a family of deformations of non-trivial bihermitian structures $(J, J^-_t, h_t)$ on $X$ by using $\b$.
In addition, if the class $[\b\cdot \ome]$ does not vanish
for a K\"ahler form $\ome$, then 
the complex structure $J_t^-$ is not equivalent to $J$ for small $t\neq 0$ under diffeomorphisms.
  Our method is based on the construction of \complex and \K\"ahler structures developed in \cite{Go1} and \cite{Go2}.
 As applications, we obtain such deformations 
 of bihermitian structures on del Pezzo surfaces, the Hirtzebruch surfaces $F_2, F_3$ and degenerate del Pezzo surfaces. 
Further we show that   
del Pezzo surfaces $S_n\, (5\leq n\leq 8)$, $F_2$ and degenerate del Pezzo surfaces admit bihermitian structures for which  $(X, J^-_t)$ is not biholomorphic to 
$(X, J)$ for small $t\neq 0$.
 \end{abstract}
\tableofcontents
\numberwithin{equation}{section}
\section*{Introduction}

A bihermitian structure on a $C^\infty$ manifold $X$ consists of a pair of integrable complex structures $J^+$ and $J^-$ with a Riemannian metric $h$ which is hermitian with respect to both $J^+$ and $J^-$. 
If a complex manifold $(X, J)$ has a bihermitian structure $(J^+, J^-, h)$ with the property $J^+=J$, then we say that $(X, J)$ admits a (compatible) bihermitian structure. 
A bihermitian structure $(J^+, J^-, h)$ is {\it  distinct} if  the complex manifold $(X, J^+)$ is not biholomorphic to
$(X, J^-)$.
We have the two $\ol\pa$-operators $\ol\pa_+$ and $\ol\pa_-$ corresponding to the complex structures $J^+$ and $J^-$ respectively. 
In this paper we always assume that a bihermitian structure satisfies the condition, 
\bgn{equation}\label{eq: torsion condition}
-d^c_+\ome_+ =d^c_-\ome_- =db,
\end{equation}
where $d^c_\pm=\sqrt{-1}(\ol\pa_\pm-\pa_\pm)$ and 
$\ome_\pm$ denote the fundamental 2-forms with respect to $J^\pm$ and 
$b$ is a real $2$-form.
(Note that if $H:=-d^c_+\ome_+ =d^c_-\ome_-$ is not $d$-exact but $d$-closed, $(J^+, J^-, h)$ is called the $H$-twisted bihermitian structure.)
There is a research of compact complex surfaces which admit bihermitian structures from the view point of Riemannian geometry
\cite{A.G.G}.
Bihermitian structures with the condition (\ref{eq: torsion condition}) appeared on the target space of $(2,2)$ supersymmetric sigma model \cite{G.H.R}. Surprisingly it turned out that
there is a one to one correspondence between \K\"ahler structures and 
bihermitian structures with the condition (\ref{eq: torsion condition})\cite{Gu1}.
It is thus expected that the construction of interesting and various \K\"ahler structures 
would be a major step of development of the theory of bihermitian structures.
Let $(X, J)$ be a compact K\"ahler manifold with a K\"ahler form $\ome$. 
In the paper \cite {Go1, Go2}, the author constructed a family of deformations of bihermitian structures by using a holomorphic Poisson structure $\b$.
In the present paper, we shall obtain another family of deformations of 
bihermitian structures $(J^+_t, J^-_t, h_t)$ of $(X, J)$, starting with the ordinary K\"ahler structure which satisfies
$J^+_t=J$ for all $t$, $J^-_0=J$ and $J^-_t\neq \pm J $ for small $t\neq 0$, 
where $t$ is a parameter of deformations.

Throughout this paper we will assume that $X$ is the underlying differential manifold of a complex manifold $M=(X, J)$ with the structure sheaf ${\cal O}_M$.
We denote by $\Theta$ the sheaf of germs of sections of the tangent bundle $T^{1,0}_J$ of $M=(X, J)$ and $\w^p\Theta$ is the sheaf of germs of $p$-th skew symmetric tensors of $\Theta$. 
Our main theorem is the following : 
\bgn{theorem}\label{th: main theorem}
Let $M=(X, J) $ be a compact K\"ahler manifold. We assume that the direct sum of cohomology groups $\oplus_{i=0}^3H^i (M, \w^{3-i}\Theta)$ vanishes. 
Then for every K\"ahler form $\ome$ and every non-zero holomorphic Poisson structure $\b$, 
there exist deformations of bihermitian structures $( J^+_t, J^-_t, h_t)$ which  satisfies,
\bgn{equation}\label{eq: type II}
J^+_t=J^-_0=J, \qquad \frac d{dt} J^-_t|_{t=0} =
-2(\b\cdot\ome+\ol\b\cdot\ome),
\end{equation}
where $\b\cdot\ome$ is the $\ol\pa$-closed forms of type $(0,1)$ with coefficients in the tangent bundle $T^{1,0}_J$ which is given by the contraction between $\b$ and $\ome$, 
and $\ol\b\cdot \ome$ is the complex  conjugate.
The $\ol\pa$-closed form $\b\cdot\ome$  gives rise to the Kodaira-Spencer class 
$-2 [\b\cdot \ome]\in H^1(M, \Theta)$
of deformations $\{J^-_t\}$.
\end{theorem}
The condition (\ref{eq: type II}) implies that $J^-_t\neq \pm J$ 
for small $t\neq 0$ as almost complex structures.
However these $J^-_t$ and $J$ might be equivalent under  diffeomorphisms.
If the class $[\b\cdot\ome]\in H^1(M, \Theta)$ does not vanish, the family of deformations $\{J^-_t\}$ is not obtained by the action of a one-parameter family of diffeomorphisms on $J$. 
Thus the complex manifold $(X, J^-_t)$ is different from $(X, J)$ for small $t\neq 0$.
\bgn{theorem}
Let $M=(X, J) $ be a compact K\"ahler manifold. 
We assume that the direct sum of cohomology groups $\oplus_{i=0}^3H^i (M, \w^{3-i}\Theta)$ vanishes and in addition, the class $[\b\cdot \ome]\in H^1(M, \Theta)$ does not vanish for a K\"ahler form $\ome$ and a holomorphic Poisson structure $\b$.
Then there exist deformations of distinct bihermitian structures $( J, J^-_t, h_t)$, that is, $(X, J^-_t)$ is not biholomorphic to $M=(X, J)$ for small $t\neq 0$.
\end{theorem}

The infinitesimal deformations of \complex structures are given by the direct sum of cohomology group,
$$
H^0(M, \w^2\Theta)\oplus H^1(M, \Theta)\oplus H^2(M, {\cal O}_X),
$$
where $H^1(M, \Theta)$ is the space of the Kodaira-Spencer classes which gives infinitesimal deformations of 
(usual) complex structures. The cohomology group $H^2(M, {\cal O}_X)$ corresponds to  the exponential action of $\ol\pa$-closed $2$-form of type $(0,2)$, which
is often called transformations by {\it $b$ fields}.
The space  $H^0(M,\w^2\Theta)$ corresponds to  deformations of \complex structures $\{ \J_{\b t}  \}$ by a Poisson structure $\b$ which are called 
{\it the Poisson deformations}. 
As in deformations of complex manifolds, there exits an obstruction to deformations of \complex structures in general.
The obstruction space to deformations of \complex structures at $\J_J$ is given by the direct sum of the ordinary cohomology groups,
$$
\bigoplus_{i=0}^3 H^i(M, \w^{3-i}\Theta):=H^0(M, \w^3\Theta)\oplus H^1(M, \w^2\Theta)\oplus H^2(M, \Theta)\oplus H^3(M,{\cal O}_M)
$$
which is the obstruction space in the theorem \ref{th: main theorem}.
If the space of the obstruction vanishes,
we can apply the method  in \cite{Go1} and \cite{Go2} 
to construct a family of \K\"ahler structures which corresponds to  the one of bihermitian structures in the theorem \ref{th: main theorem}.
More precisely,
the complex structure $J$ gives a \complex structure $\J_J$ and the K\"ahler  structure $\ome$ also provides
the $d$-closed non-degenerate, pure spinor $\psi=e^{\sqrt{-1}\ome}$
which induces the \complex structure $\J_{\psi}$.  The pair $(\J_J, \psi)$ gives rise to a \K\"ahler structure
$(\J_{J}, \J_{\psi})$.
It is the essential feature that 
the generalized geometry inherits the symmetry of the Clifford group of the direct sum of the tangent bundle $T$ and the cotangent bundle $T^*$ on a manifold $X$.
The space of almost \K\"ahler structures forms 
an orbit by the diagonal action of the Clifford group. 
Thus we construct deformations of almost \K\"ahler structures with one pure spinor staring with $(\J_J, \psi)$ by the action of the Clifford group, 
$$
\J_t= \Ad_{e^{Z(t)}}\J_J, \quad \psi_t= e^{Z(t)}\psi,
$$  
where $e^{Z(t)}$ is a family of the Clifford group  and $\Ad_{e^{Z(t)}}$ denotes the adjoint action of the Clifford group on $\J_J$ (see \cite{Go1}).
The pair $(\J_t, \psi_t)$ induces the almost \K\"ahler structure
$(\J_t, \J_{\psi_t})$ and
then the corresponding bihermitian structures $(J^+_t, J^-_t, h_t)$ are given by the action of $\Gam^\pm_t\in$ GL$(TX)$ by
$J^\pm_t= (\Gam^\pm_t)^{-1}\circ J\circ \Gam^\pm_t$, where 
$\Gam^\pm_t$ is explicitly described in terms of $Z(t)$ (see \ref{eq: Gam t pm } in section 3).
Thus our problem is reduced to construct $Z(t)$ which satisfies the following three conditions:
\bgn{align}
&\J_t:=\Ad_{e^{Z(t)}}\J_ J\,\, \text{\rm are integrable \complex structures}\\
&d\psi_t:=d e^{Z(t)}\psi =0 
\\
& (\Gam^+_t)^{-1}\circ J\circ \Gam^+_t= J
\end{align}
Then $Z(t)$ yields deformations $(\J_t, \J_{\psi_t})$ of \K\"ahler structures which gives rise to 
bihermitian structures in the theorem \ref{th: main theorem} 
(see section 4 for more detail). \footnote{Note that if $\psi_t$ is closed, the induced structure $\J_{\psi_t}$ is integrable.}

Hitchin \cite{Hi3} constructed deformations of bihermitian structure of the type $(J, J^-_t, h_t)$ by the Hamiltonian diffeomorphisms 
on del Pezzo surfaces and Gualtieri \cite{Gu3} extended the approach to higher dimensional Poisson manifolds.
The bihermitian structures which they constructed give the equivalent two complex structures under diffeomorphisms. Our constructions enable us to obtain distinct bihermitian structures.
\footnote{\noindent \noindent The author received a note that Gualtieri also developed a modified approach to obtain bihermitian structures on $F_2$ recently.}

In section 1, we will give a short explanation of deformations of \complex structures. 
 Deformations of \complex structures are often described in the language of 
complex Lie algebroid \cite{L.W.X}, \cite{Gu1}. It is necessary to translate it in the terms of the action of the (real) Clifford group for our construction of \K\"ahler structures.
In section 2, we recall the stability theorem of \K\"ahler structure with one pure spinor which was shown in \cite{Go1} and \cite{Go2}.
In section 3, we give a description of $\Gam_t^\pm$ which gives deformations of bihermitian structures corresponding to the ones of \K\"ahler structures. 
In section 4, we will construct deformations of bihermitian structures in the main theorem \ref{th: main theorem} as formal power series. 
In section 5, we will show the convergence of the power series constructed in section 4  and finish our proof of the main theorem. 
In section 6, we apply our method to complex surfaces.
In the case of complex surfaces, we only need to show that the cohomology groups
$H^1(M, K_M^{-1})$ and $H^2(M, \Theta)$ vanish
to obtain deformations in the theorem \ref {th: main theorem}, where $K_M$ is the canonical line bundle.

In subsection 6.1, we show that every del Pezzo surface admits deformations of bihermitian structures as in theorem \ref{th: main theorem}.
 Let $S_n$ be a del Pezzo surface which is the blow-up of $\C P^2$ at $n$ points. Then we prove that if $n\geq 5$, there exists a class $[\b\cdot \ome]\in H^1(S_n, \Theta)$ which does not vanish for a K\"ahler form $\ome$.  As a result, we obtain distinct bihermitian structures on $S_n$ $(n\geq 5)$.
 In subsection 6.2, we will give several vanishing theorems of $H^1(M, K_M^{-1})$ and $H^2(M, \Theta)$ on a complex surface $M$.
 In subsection 6.3 we will show the non-vanishing theorem of the class $[\b\cdot \ome]\in H^1(M, \Theta)$ which gives rise to unobstructed deformations. 
 Applying these vanishing theorems and the non-vanishing theorem, we obtain 
bihermitian structure $(J^+, J^-)$ on $F_2=(X, J)$ on which the complex manifold $(X, J^+)$ is $F_2$ and 
$(X, J^-)$ is $\C P^1\times \C P^1$ in subsection 6.4.
We also show that the Hirtzebruch surface $F_3$ admits bihermitian structures.
{\it Degenerate del Pezzo surfaces} are the blow-up of $\C P^2 $ at $r$ points, $0\leq r \leq 8$ which are in {\it almost general position}
(see subsection 6.5 for more detail). It turns out that the obstruction spaces still vanish on degenerate del Pezzo surfaces.
If the anti-canonical line bundle is not ample, then there is a $(-2)$-curve $C$ with $K\cdot C=0$
and it follows that the  class $[\b\cdot\ome]$
does not vanish. 
Hence we obtain bihermitian structures in theorem \ref{th: main theorem} on the degenerate del Pezzo surfaces which yield distinct two complex manifolds.
We contract all $(-2)$-curves on a degenerate del Pezzo to obtain a del Pezzo surface with rational double points, which is called the {\it Gorenstein log del Pezzo surface}, \cite{DV, De, AN}.
 In appendix I, we give the power series construction of the Kuranishi family of \complex structures. 
In appendix II, we collect necessary formulae and give an explanation of the Schouten bracket and the Jacobi identity of the brackets.

The author would like to thank for Professor A. Fujiki for valuable comments. 
He is also thankful to Professor V. Apostolov and Professor M. Gualtieri for our remarkable discussion at Montreal.

\section{Deformations of \complex structures}
Let $\J$ be a \complex structure on a compact manifold $X$ of real dimension $2n$. 
Then the \complex structure $\J$ gives the decomposition 
$(\TT)^\C=L_\J\oplus\ol{L}_\J$, where $L_\J$ denotes the eigenspace with eigenvalue $\sqrt{-1}$ and $\ol{L}_\J$ is the complex conjugate of $L_\J$. 
For a section $\e\in\w^2\ol{L}_\J$, the exponential $e^{\e}$ is regarded as the section of 
complex Clifford group which also induces the section of SO$(\TT, \C)$ by the adjoint action Ad$_{e^\e}$.
Small deformations of almost \complex structure $\J$ are written by a section of $\w^2\ol{L}_\J$, 
$$\J_t=\Ad_{e^{\e(t)}}J,$$where $t$ is the parameter of deformations.
The integrability condition of almost \complex structure $\J_t$ is given by 
the Maurer-Cartan equation, 
\bgn{equation}\label{eq: Maurer-Cartan }
d_L\e(t)+\frac12[\e(t), \e(t)]_S=0,
\end{equation}
where $d_L : \w^k\ol{L}_\J\to \w^{k+1}\ol{L}_\J$ denotes the exterior derivative of the Lie algebroid $L_\J$ and the bracket $[\,,\,]_S$ is the Schouten bracket of $\ol{L}_\J$.
Hence the problem of deformations reduces to solving the Maurer-Cartan equation, 

At first, we write a family of sections $\e(t)$ as a power series in $t$
\bgn{equation}\label{eq: e(t), power series}
\e(t)=\e_1t+\e_2\frac {t^2}{2!}+\cdots,
\end{equation}
Note that our power series starts from $\e_1$. 
Substituting the power series (\ref{eq: e(t), power series})
into the \MC, we obtain the equation on $t$.  
We denote by $\([\e(t), \e(t)]_S\)_{[k]}$ the $k$ th homogeneous term in $t$. Then the equation is reduced to infinitely many equations, 
\bgn{equation}\label{eq: k th term MC}
\frac 1{k!}d_L\e_k +\frac12\([\e(t), \e(t)]_S\)_{[k]}=0
\end{equation}
We have the differential complex $(\w^{\bullet}\ol L_\J, d_L)$ which is elliptic, 
$$
0\to \ol L_\J \overset{d_L}\to \w^2\ol L_\J\overset{d_L}\to
\w^3\ol L_\J\overset{d_L}\to \cdots
$$
For $k=1$, the equation is $d_L\e_1=0$. 
It implies that $\e_1$ is a section of $\w^2 \ol L_\J$ which is $d_L$-closed.
We take $\e_1$ as a Hormonic section which satisfies 
$$
\trian_L\e_1=(d_Ld_L^*+d_L^*d_L)\e_1=0,
$$
where $d_L^*$ is the formal adjoint operator of $d_L$ with respect to a
Riemannian metric on $M$.
There are actions of diffeomorphisms and $d$-exact $b$-fields on \complex structures
which generate $d_L$-exact sections of $\w^2 L_\J$ infinitesimally.
We identify deformations by both actions of diffeomorphisms and $d$-exact $b$-fields.
It implies that the infinitesimal deformations (the first order deformations) are given by the cohomology group $H^2(\ol L_\J)$ of the elliptic differential complex $(\w^{\bullet}\ol L_\J, d_L)$ . 
The third cohomology group $H^3(\ol L_\J)$ is regarded as the space of the obstructions to deformations. 
The deformation theory of \complex structures was already discussed in \cite{Gu1} by the implicit function theorem.
We will give the different construction of deformations of \complex structures by using the power series, which is analogous to the one of original Kodaira-Spencer theory.
Our method yields an estimate of the convergent series which is necessary for the construction of \K\"ahler and bihermitian structures.

\bgn{theorem}\label{th: deform of G.C}
If the cohomology group $H^3(\ol L_\J)$ vanishes, then we have 
a family of deformations of \complex structures which are parametrized by an 
open set of $H^2(\ol L_\J)$.
\end{theorem}
\bgn{proof}
We solve the equations (\ref{eq: k th term MC}) by the induction on the degree of $t$.
We assume that there are sections $\e_1,\cdots, \e_{k-1}\in\w^2\ol L_\J$ which satisfy the equations 
\bgn{equation}\label{eq: degree i equation}
\frac1{i!}\e_{i}+\frac12\( [\e(t), \e(t)]_S\)_{[i]}=0,
\end{equation}
for all $i<k$. 

Then we shall show that there is a section $\e_k$ which satisfies 
the equation (\ref{eq: k th term MC}). 
 The Schouten bracket and the Lie algebroid derivative $d_L$ satisfy the following relations for sections $\e_1, \e_2\in \w^\bullet\ol L_\J$,
 \bgn{proposition}\label{prop: Schouten relations}
 \bgn{align*}
 &[\e_1, \e_2]_S=(-1)^{|\e_1|\,|\e_2|} [\e_2, \e_1]\\ \\
 &d_L[ \e_1, \e_2]_S=[d_L\e_1, \e_2]_S+(-1)^{|\e_1|}[\e_1, d_L\e_2]_S\\ \\
 &(-1)^{|\e_1|\,|\e_3|}\[[\e_1, \e_2],\,\e_3\]
 +(-1)^{|\e_2|\,|\e_1|}\[[\e_2, \e_3],\, \e_1\]
 +(-1)^{|\e_3|\,|\e_1|}\[[\e_3, \e_1],\,\e_2\]=0
 \end{align*}
 where we denote by $|\e_i|$ the degree of $\e_i$.
 \end{proposition}
 The relations in the proposition \ref{prop: Schouten relations} are already known 
\cite{KoSc}.
 Note that in our case the degree of $\e_i$ is even and we have the ordinary 
 Jacobi identity.
 The $k$ th order term of the Schouten bracket $[\e(t), \e(t)]_S$ is given by 
 \bgn{align}\label{ali: degree k Schouten}
\( [\e(t), \e(t)]_S\)_{[k]}=&\sum_{\stackrel {\scriptstyle i+j=k}{0<i, \,j<k}}\frac 1{i!}\frac 1{j!}
\[\e_{i},\, \e_{j}\]_S.
 \end{align}
 Then substituting (\ref{eq: degree i equation}) into (\ref{ali: degree k Schouten}) and applying the proposition \ref{prop: Schouten relations}, we have 
 \bgn{align}
d_L \( [\e(t), \e(t)]_S\)_{[k]}=&\sum_{\stackrel {\scriptstyle i+j=k}{0<i, \,j<k}}\frac1{i!}\frac1{j!}
d_L\[\e_i,\, \e_{j}\]_S\\
=&2\sum_{\stackrel {\scriptstyle i+j=k}{0<i, \,j<k}}\frac1{i!}\frac1{j!}\[d_L\e_{i}, \, \e_{j}\]_S\\
=&-\sum_{\stackrel {\scriptstyle l+m+j=k}{0<l, \,m, \,j<k}}\frac1{l!}\frac1{m!}\frac1{j!}
\[[\e_l, \e_m]_S,\, \e_j\]_S=0
 \end{align}

Hence $\( [\e(t), \e(t)]_S\)_{[k]}\in\w^2\ol L_\J$ is $d_L$-closed which is a representative of the cohomology class in $H^3(\ol L_\J)$.
Since we assume that $H^3(\ol L_\J)$ vanishes, we have a solution $\e_k$ of the equation (\ref{eq: k th term MC}).
We use a Riemannian metric on $X$ to construct $\e_k$ by using the formal adjoint operator $d_L^*$ and the Green operator $G_L$ of the elliptic complex $(\w^\bullet\ol L_\J, d_L)$
\bgn{equation}\label{eq: Greeen }
\frac1{k!}\e_{k}=-\frac12d_L^*G_L\( [\e(t), \e(t)]_S\)_{[k]}
\end{equation}
In fact, it follows from the Hodge decomposition theorem that $\e_{k}$ satisfies the equation (\ref{eq: k th term MC}).
Thus we obtain the solution $\e(t)$ of the {\MC} as a formal power series.
In order to show the power series $\e(t)$ is a convergent series which is further 
a smooth solution, we apply the standard method due to Kodaira-Spencer. 
Let $P(t)=\sum_k P_k t^k$ be a power series in $t$ whose coefficients are
sections of a vector bundle with a metric.
We denote by $\|P_k\|_s$ the Sobolev norm of the section $P_k$ 
which is given by the sum of the $L^2$-norms of  $i$ th derivative of $P_k$ for all $i\leq s$, where $s$ is a positive integer with $s>2n+1$.
We put $\|P(t)\|_s=\sum_k \|P_k\|_st^k$.
Given two power series $P(t), Q(t)$, if 
$\|P_k\|\leq \|Q_k\|$ for all $k$, then we denote it by 
$$
P(t)\<Q(t).
$$
For a positive integer $k$, if $\|P_i\|\leq \|Q_i\|$ for all $i\leq k$,
we write it by 
$$
P(t)\underset k{\<} Q(t).
$$
We also use the following notation.
If $P_i=Q_i$ for all $i\leq k$, we write it by 
\bgn{equation}\label{eq: notation 1}
P(t)\underset k{\equiv} Q(t).
\end{equation}
Let $M(t)$ be a convergent power series defined by 
\bgn{equation}\label{eq: M(t)}
M(t) =\sum_{\nu=1}^\infty \frac1{16c}\frac{(ct)^\nu}{\nu^2}=\sum_{\nu=1}^\infty M_\nu t^\nu,
\end{equation}
for a positive constant $c$, which is determined later suitably.
The key point is the following inequality,
\bgn{equation}\label{eq: M(t)}
M(t)^2 \< \frac 1c M(t)
\end{equation}
We put $\lam =c^{-1}$. 
Then we also have 
\bgn{equation}\label{e M(t)}
e^{M(t)}\<\frac 1\lam e^\lam M(t).
\end{equation}
We assume that our power series $\e(t)$ satisfies the inequality for an integer $k>1$, 
\bgn{equation}
\|\e(t)\|_s\underset{k-1}\< M(t)
\end{equation}
We apply the standard estimate of elliptic differential operators to obtain an estimate of  
the solution $\e_{[k]}$ in (\ref{eq: Greeen }), 
\bgn{align}
\frac2{k!}\|\e_k\|_s\leq& C_1\|\(\[\e(t), \,\e(t)\]_S\)_{[k]}\|_{s-1}
=\sum_{\stackrel {\scriptstyle i+j=k}{0<i, \,j<k}}\frac1{i!}\frac1{j!}
C_1\|\[\e_{i},\, \e_{j}\]_S\|_{s-1}\\
\leq &2C_1\sum_{\stackrel {\scriptstyle i+j=k}{0<i, \,j<k}}\frac1{i!}\frac1{j!}\|\e_i\|_s\,
\|\e_j\|_s=2C_1\sum_{\stackrel {\scriptstyle i+j=k}{0<i, \,j<k}}
M_iM_j\\
\leq &2C_1\lam M_k
\end{align} 
Hence if we choose a constant $\lam$ with $C_1\lam <1$, 
then it follows that $\frac1{k!}\|\e_k\|_s< M_k$. 
It implies that $\e(t)$ is a convergent series. 
From our construction, the series $\e(t)$ satisfies 
\bgn{equation}\label{eq: e(t)}
\e(t) =\e_1t-\frac12 d_L^*G_L[\e(t), \e(t)]_S
\end{equation}
Since $H^3(\w^\bullet\ol L_\J)=\{0\}$, we have a differential equation 
$$
\trian_L\e(t) +\frac12 d_L^*[\e(t), \e(t)]_S=0,
$$
which is elliptic for sufficiently small $\e(t)$. 
Thus it follows that $\e(t)$ is a smooth solution.
\end{proof}
In the appendix, we further construct the Kuranishi family of deformations of \complex structures which gives the space of deformations even in the cases where 
the obstruction space $H^3(\w^\bullet\ol L_\J)$ does not vanish.

We denote by CL the Clifford algebra bundle of $\TT$ on a manifold $X$ which admits filtrations of even degree and odd degree,
\bgn{align*}
&\CL^1\subset \CL^3\subset\CL^5\subset\cdots\\
&\CL^0\subset \CL^2\subset\CL^4\subset\cdots,
\end{align*}
where $\CL^0=\TT$ and $\CL^2$ denotes the the subbundle of CL which consists 
of elements of degree $2$ or $0$,
(for simplicity, we call CL the Clifford algebra of $\TT$.)
Let $\J$ be a \complex structure on $X$ which gives the decomposition, 
$$
(\TT)^{\C}=L_\J\oplus\ol{L}_\J
$$
We denote by $\w^p\ol L_\J$ the bundle of the $p$ th skew symmetric tensor of $\ol L_\J$. 
Let $U^{-n}_\J$ be the line bundle of $(X,\J)$ which consists of non-degenerate, complex pure spinors corresponding to $\J$.
We call $U_\J^{-n}$ the canonical line bundle $K_\J$ of $\J$.
There is the action of $\TT$ on differential forms $\w^\bullet T^*$ by the interior product and the exterior product which induces the spin representation of 
the Clifford algebra CL on $\w^\bullet T^*$.
By the action of $\w^p\ol L_\J$ on $K_\J$, we have the vector bundles, 
$$
U^{-n+p} :=\w^p\ol L_\J\cdot K_\J
$$
Then the differential forms $\w^\bullet T^*$ on $X$ are decomposed into 
$$
\w^\bullet T^*=\bigoplus_{p=0}^{2n} U^{-n+p}
$$
We denote by $\pi_{U^{-n+p}}$ the projection to the bundle $U^
{-n+p}$.
The set of almost \complex structures forms the orbit of
the (real) Clifford group of the Clifford algebra CL of $\TT$ which acts on $\J$ by 
the adjoint action.
The Lie algebra of the Clifford group is the subalgebra CL$^2$. 
Small deformations of almost \complex structures $\{\J_t\}$ are given in terms of the
adjoint action,
$$
\J_t:=\Ad_{e^{a(t)}}\J,
$$
where $a(t)=a_1t+\frac 1{2!}a_2t^2+\cdots$ is a CL$^2$-valued power series in $t$.

In order to obtain deformations of \K\"ahler structures, we need to consider 
a section $a(t)$ of the bundle CL$^2(\TT)$. 

it is crucial that
the set of almost \K\"ahler structures just forms an orbit of the action of 
the real Clifford group and 
deformations of almost \K\"ahler structures are not given by the action of 
the complex Clifford group.
The following lemma is necessary for the construction of \K\"ahler structures.
which is already proved in \cite{Go1}.
\bgn{lemma}\label{lem: e(t)-> a(t)}
For small deformations of almost \complex structures given by 
$\J_t:=\Ad_{e^{\e(t)}}\J_0$ as before, there exists a unique family of sections $a(t)$ of 
real Clifford bundle CL$^2$ such that
$$
\J_t=\Ad_{e^{a(t)}}\J_0,
$$
and $a(t) $ is in the real part of $\w^2\ol L_\J\oplus \w^2 L_\J$.
Conversely, If we have a family of deformations of almost \complex structure $\J_t
=\Ad_{e^{a(t)}}\J_0$ which is given by the action of a section $a(t) \in \CL^2$, 
then there exists a unique section $\e(t)\in \w^2\ol L_\J$ such that 
$\J_t=\Ad_{e^{\e(t)}}\J_0$.
\end{lemma}
We consider the operator $e^{-a(t)}\circ d\circ e^{a(t)}$ acting on $K_J=U^{-n}$.
Then as discussed in \cite{Go0}, the operator $e^{-a(t)}\circ d\circ e^{a(t)}$
is a Clifford-Lie operator of order $3$ whose image is in $U^{-n+1}\oplus U^{-n+3}$.

It is shown in \cite{Go1} that the almost \complex structure $\J_t=\Ad_{e^{a(t)}}\J$ is integrable if and only if 
the projection to the component $U^{-n+3}$ vanishes, that is,
$$
\pi_{U^{-n+3}}e^{-a(t)}\circ d\circ e^{a(t)}=0
$$
We denote by $\( \pi_{U^{-n+3}}e^{-a(t)}\circ d\circ e^{a(t)}\)_{[k]}$ the $k$ th term of $\pi_{U^{-n+3}}e^{-a(t)}\circ d\circ e^{a(t)}$.
Let $\J_J$ be the \complex structure on $X$ defined by a ordinary complex structure $J$. We put $M=(X, J)$.
Then as in \cite{Go1} and \cite{Go2}, the obstruction space to  deformations of $\J_J$ is given by $\oplus_{p+q=3}H^p(M, \w^q\Theta)$.
In this case, the canonical line bundle is the ordinary one $K_J$ which consists of complex forms of type $(n,0)$.
Thus by the theorem \ref{th: deform of G.C} and the lemma \ref{lem: e(t)-> a(t)}, we have the following,

\bgn{proposition}\label{prop: k th-G.complex}
Let $M=(X, J)$ be a compact K\"ahler manifold with a K\"ahler form $\ome$.
We assume that the cohomology groups $\oplus_{p+q=3}H^p(M, \w^q\Theta)
$ vanish. 
If there is a set of sections $a_1,\cdots a_{k-1}$ of CL$^2$ which satisfies 
\bgn{equation}\label{eq: k th-G.complex}
\(\pi_{U^{-n+3}}e^{-a(t)}\, de^{a(t)}\)_{[i]}=0, \qquad \text{\rm for all }
i<k,
\end{equation}
and $\|a(t)\|_s\underset{k-1} \<C_1M(t)$,
then there is a section $a_k$ of CL$^2$ which satisfies the followings:
$$
\pi_{U^{-n+3}}\(e^{-a(t)}\, de^{a(t)}\)_{[k]}=0
$$
and $\|a(t)\|_s\underset{k}\< C_1M(t)$, where 
$a(t) =\sum_{i=1}^\infty \frac 1{i!} a_i t^i$ and 
$M(t)$ is the convergent series in (\ref{eq: M(t)}) and
$C_1$ is a positive constant.
\end{proposition}
\bgn{proof}
We use the notation as in (\ref{eq: notation 1}). 
The equation (\ref{eq: k th-G.complex}) is equivalent to say that 
there is a section $\h E(t) \in \TT$ such that , 
\bgn{equation}\label{proof: k th-G.complex-1}
e^{-a(t)}\, de^{a(t)}\cdot\phi \underset{k-1}{\equiv} \h E(t)\cdot\phi,
\end{equation}
for all $\phi\in K_J$.
By the left action of $e^{a(t)}$ on both sides of the equation (\ref{proof: k th-G.complex-1}), we have 
$$
de^{a(t)}\cdot\phi\underset{k-1}{\equiv} e^{a(t)}\h E(t)\cdot\phi.
$$
We put $E(t) =e^{a(t)}\h E(t) e^{-a(t)}$. Then it follows that 
\bgn{equation}\label{proof: k th-G.complex-2}
de^{a(t)}\cdot\phi\underset{k-1}{\equiv}  E(t)\cdot e^{a(t)}\phi.
\end{equation}
From the lemma \ref{lem: e(t)-> a(t)}, we have $\e(t)\in \w^2\ol L_\J$ such that 
\bgn{equation}\label{proof: k th-G.complex-3}
e^{\e(t)}\cdot\phi\underset{k-1}{\equiv} e^{a(t)}\cdot\phi.
\end{equation}
Substituting (\ref{proof: k th-G.complex-3}) into (\ref{proof: k th-G.complex-2}), we obtain 
\bgn{equation}\label{proof: k th-G.complex-4}
de^{\e(t)}\cdot\phi\underset{k-1}{\equiv}  E(t)\cdot e^{\e(t)}\phi.
\end{equation}
By the right action of $e^{-\e(t)}$ on (\ref{proof: k th-G.complex-4}) again, 
we have 
\bgn{equation}\label{proof: k th-G.complex-5}
e^{-\e(t)}de^{\e(t)}\cdot\phi\underset{k-1}{\equiv}  
e^{-\e(t)}E(t)\cdot e^{\e(t)}\phi \underset{k-1}{\equiv}  \til E(t)\cdot\phi,
\end{equation} 
where $\til E(t) =e^{-\e(t)}E(t)\cdot e^{\e(t)}$.
Thus as in \cite{Go1}, the equation(\ref{proof: k th-G.complex-5}) 
is equivalent to the \MC ,
$$
d_L\e(t)+\frac12[\e(t), \e(t)]_S \underset{k-1}{\equiv}0.
$$
Then as is shown in the theorem \ref{th: deform of G.C}, there is a section $\e_k$ such that 
$$
d_L\e(t)+\frac12[\e(t), \e(t)]_S \underset{k}{\equiv}0.
$$
We define $a_k$ by $a_k=\e_k+\ol\e_k$. 
Then it follows that 
$$
e^{-a(t)}\, de^{a(t)}\cdot\phi \underset{k}{\equiv} \h E(t)\cdot\phi.
$$
Hence we have $\(\pi_{U^{-n+3}}e^{-a(t)}\, de^{a(t)}\)_{[k]}=0$ and 
$\|a(t)\|_s\underset{k}\< C_1M(t)$.
\end{proof}

\section{Deformations of generalized K\"ahler structures}
Let $(X, J,\ome)$ be a compact K\"ahler manifold
and $(\J, \J_{\psi})$ the \K\"ahler structure induced from
$(J,\ome)$ by $\J=\J_J$ and $\psi=e^{\sqrt{-1}\ome}$.
Since two \complex structures $\J$ and $ \J_{\psi}$ are commutative, 
the \K\"ahler structure $(\J, \J_{\psi})$ gives the simultaneous  decomposition of 
$(\TT)^\C$, 
$$
(\TT)\C=L^+_\J\oplus L^-_\J\oplus \ol L^+_\J\oplus\ol L^-_\J,
$$
where $L^+_\J\oplus L^-_\J$ is the eigenspace with eigenvalue $\sqrt{-1}$ with respect to $\J$ and $L^+_\J\oplus \ol L^-_\J$ is the eigenspace with eigenvalue $\sqrt{-1}$ with respect to $\J_\psi$ and $\ol L^\pm_\J$ denotes the complex conjugate.
In \cite{Go1,Go2}, the author showed the stability theorem of 
\K\"ahler structures with one pure spinor, which implies that 
if there is a one dimensional  analytic deformations of \complex structures $\{\J_t\}$ parametrized by $t$,
then there exists a family of non-degenerate, $d$-closed pure spinor $\psi_t$
such that the family of pairs $(\J_t, \psi_t)$ becomes deformations of 
{\K\"ahler} structures starting from $(\J, \psi)=(\J_0, \psi_0)$.
As in section 2, small deformations $\J_t$ can be written by the adjoint action of $a(t)$ in CL$^2$, 
$$
\J_t:=\Ad_{e^{a(t)}}\J_0. 
$$
Then we can obtain a family of real sections $b(t)$ of the bundle
$(L^-_{\J_0}\cdot\ol L^+_{\J_0}\oplus \ol L^-_{\J_0}\cdot L^+_{\J_0})$ such that $\psi_t=e^{a(t)}e^{b(t)}\psi_0$ is the family of non-degenerate, $d$-closed pure spinor $\psi_t$.  
The bundle $K^1=U^{0,-n+2}$  is generated by the action of real sections of $(L^-_{\J_0}\cdot\ol L^+_{\J_0}\oplus \ol L^-_{\J_0}\cdot L^+_{\J_0})$ on $\psi$ (see page 125 in \cite{Go2} for more detail).

We define $Z(t)$ by 
$$
e^{Z(t)}= e^{a(t)}\, e^{b(t)}.
$$
Since $\Ad_{e^b(t)}\J_0=\J_0$, we obtain 
$\J_t=\Ad_{e^{a(t)}}\J_0=\Ad_{e^{a(t)}}\, \Ad_{e^{b(t)}}\J_0=\Ad_{e^{Z(t)}}\J_0 $. Then the family of deformations of \K\"ahler structures is given by the action of $e^{Z(t)}$, 
$$
(\J_t, \,\psi_t) =\(\Ad_{e^{Z(t)}}\J_0, \, \,\,e^{Z(t)}\cdot\psi \,\).
$$
By the similar method as in \cite{Go2} together with the proposition \ref{prop: k th-G.complex}, we obtain the following proposition, 
\bgn{proposition}\label{prop: GK deformations}
Let $(X, J,\ome)$ be a compact K\"ahler manifold. 
We assume that the cohomology groups $\oplus_{p+q=3}H^p(X, \w^q\Theta)$ vanish. 
If there is a set of sections $a_1,\cdots a_{k-1}$ of CL$^2$ which satisfies 
$$
\pi_{U^{-n+3}}\(e^{-a(t)}\, de^{a(t)}\)_{[i]}=0, \qquad \text{\rm for all }
i<k,
$$
and $\|a(t)\|_s \underset{k-1}\< K_1M(t)$ for a positive constant $K_1$,
then there is a set of real sections $b_1,\cdots, 
b_k$ of the bundle 
$(L^-_{\J_0}\cdot\ol L^+_{\J_0}\oplus \ol L^-_{\J_0}\cdot L^+_{\J_0})$ which satisfies 
the following equations:
\bgn{align}
&\pi_{U^{-n+3}}\(e^{-Z(t)}\, de^{Z(t)}\)_{[k]}=0\\
&\( d e^{Z(t)}\cdot\psi_0\)_{[i]}=0,\qquad \text{\rm for all } i\leq k\\
&\|a(t)\|_s\underset {k}\< K_1\lam M(t)\\
&\|b(t)\|_s\underset{k}\<K_2M(t)
\end{align}
where $a_k$ is the section constructed in the proposition \ref{prop: k th-G.complex} and
$e^{Z(t)}=e^{a(t)}\, e^{b(t)}$ and $M(t)$ is the convergent series in (\ref{eq: M(t)}) and 
a positive constant $K_2$ is determined by $\lam$ and $K_1$.
The constant $\lam$ in $M(t)$ will be suitably selected to show the convergence of the power series $Z(t)$ in section 6.
\end{proposition}
\section{Deformations of bihermitian structures}
 We use the same notation as in pervious sections.
 There is a one to one correspondence between \K\"ahler structures and bihermitian structures with 
 the condition (\ref{eq: torsion condition}).
 In this section we shall give an explicit description of $\Gam^\pm_t$ which givens rise to
 bihermitian structure $(J^+_t, J^-_t)$ 
 corresponding to deformations $(\J_t, \psi_t)$.
 The correspondence is defined at each point on a manifold, that is, the correspondence between tensor fields
 which allows us to obtain almost bihermitian structures from almost \K\"ahler structures.
 The non-degenerate, pure spinor $\psi_t$ induces the \complex structure $\J_{\psi_t}$. 
 Since $(\J_t, \J_{\psi_t})$ is a \K\"ahler structure and $\J_t$ commutes with $\J_{\psi_t}$, 
 we have the simultaneous decomposition of $(\TT)^\C$ into four eigenspaces,  
 $$
 (\TT)^\C=L^+_{\J_t}\oplus L^-_{\J_t}\oplus \ol {L_{\J_t}^+}\oplus \ol {L_{\J_t}^-},
 $$
where each eigenspace is given by the intersection of eigenspaces of both $\J_t$ and $\J_{\psi_t}$, 
 \bgn{align*}
 &L_{\J_t}^-=L_{\J_t}\cap \ol L_{\psi_t},\qquad \ol {L_{\J_t}^+}=\ol L_{\J_t}\cap \ol L_{\psi_t}\\
 &L_{\J_t}^+=L_{\J_t}\cap L_{\psi_t},\qquad\ol {L_{\J_t}^-}=\ol L_{\J_t}\cap L_{\psi_t},
   \end{align*}
where $ L_{\J_t}$ is the eigenspace of $\J_t$ with eigenvalue $\sqrt{-1}$
and $L_{\psi_t}$ denotes the eigenspace of $\J_{\psi_t}$ with eigenvalue $\sqrt{-1}$.
Since $\J_t=\Ad_{e^{Z(t)}}\J_0=e^{Z(t)}\, \J_0\, e^{-Z(t)}$ and 
$\J_{\psi_t}=\Ad_{e^{Z(t)}}\J_\ome$, we have the isomorphism between eigenspaces,
  $$
  \Ad_{e^{Z(t)}} : \ol{L_{\J_0}^\pm}\to \ol{L_{\J_t}^\pm}.
  $$
Let $\pi$ be the projection from $\TT$ to the tangent bundle $T$. We restrict the map $\pi$ to the 
eigenspace $\ol{L_{\J_0}^\pm}$ which yields the map $\pi_t^\pm : \ol{L_{\J_t}^\pm}\to T^\C$.
Let $T^{1,0}_{J^\pm_t}$ be the complex tangent space of type $(1,0)$ with respect to 
  $J^\pm_t$. 
  Then it follows that $T^{1,0}_{J^\pm_t}$ is given by the image of $\pi_t^\pm$, 
  $$
  T^{1,0}_{J^\pm_t}=\pi_t^\pm( \ol{L^\pm_{\J_t}})
  $$
Since deformations of \K\"ahler structures are given by the action of $e^{Z(t)}$, 
the ones of bihermitian structures $J^\pm_t$ should be described by the action of $\Gam^\pm_t\in $GL$(T)$ 
which is obtained from $Z(t)$. 
We shall describe $\Gam^\pm_t$ in terms of $a(t)$ and $b(t)$.
A local basis of $\ol{L_{\J_0}^\pm}$ is given by 
$$\{ \Ad_{e^{\pm\sqrt{-1}\ome}}V_i= V_i\pm \sqrt{-1}[\ome, V_i]\, \}_{i=1}^n,$$
for a local basis $\{ V_i\}_{i=1}^n$ of $T^{1,0}_J$, 
where we regard $\ome$ as an element of the Clifford algebra and then the bracket 
$[\ome, V_i]$ coincides with the interior product $i_{V_i}\ome$.
It follows that the inverse map $(\pi^\pm_0)^{-1} : T^{1,0}_{J} \to \ol {L_{\J_0}^\pm}$ 
is given by the adjoint action of $e^{\pm\sqrt{-1}\ome}$, 
\bgn{equation}\label{eq: pi0{-1}}
\Ad_{e^{\pm\sqrt{-1}\ome}}=(\pi^\pm_0)^{-1}.
\end{equation}
 We define a map $(\Gam_t^\pm)^{1,0}: T^{1,0}_J \to T^{1,0}_{J^\pm_t}$ by the composition, 
 \bgn{align}
 (\Gam_t^\pm)^{1,0}=&\pi^\pm_t\circ \Ad_{e^{Z(t)}}\circ (\pi^\pm_0)^{-1}\\
 =&\pi\circ \Ad_{e^{Z(t)}}\circ \Ad_{e^{\pm\sqrt{-1}\ome}}
\end{align}
\medskip
$$\xymatrix{
\ol{L^\pm_{{\cal J}_0}}\ar@{->}[rr]^{\text{\rm Ad}_{e^{Z(t)}} }\ar@{->}[d]_{\pi^\pm_0}& &\ol{L^\pm_{{\cal J}_t}}\ar@{->}[d]^{\pi^\pm_t}\\
T_J^{1,0}\ar@{->}[rr]_{(\Gamma^\pm_t)^{1,0}}&& T^{1,0}_{J_t^\pm}
}
$$

Together with the complex conjugate $(\Gam^\pm_t)^{0,1} :T^{0,1}_J \to T^{0,1}_{J^\pm_t}$, 
we obtain the map $\Gam^\pm_t$ which satisfies 
$J^\pm_t= (\Gam_t^\pm)^{-1} \circ J\circ \Gam_t^\pm$. 

Let $J^*$ be the complex structure on the cotangent space $T^*$ which is given by 
$\lan J^*\eta,  v\ran=\lan \eta, Jv\ran$, where $\eta\in T^*$ and $v\in T$ and 
$\lan \,,\,\ran$ denote the coupling between $T$ and $T^*$.
We define a map $\h J^\pm : \TT \to \TT$ by $\h J^\pm
(v,\eta) = v\mp J^*\eta$ for $v\in T$ and $\eta\in T^*$.
Then $\Gam^\pm_t$ is written as 
\bgn{equation}\label{eq: Gam t pm }
\Gam^\pm_t= \pi\circ \Ad_{e^{Z(t)}}\circ \h J^\pm\circ\Ad_{e^\ome}.
\end{equation}
The $k$ th term of $\Gam^\pm_t$ is denoted by $(\Gam_t^\pm)_{[k]}$ as before. 
Note that $(\Gam_t^\pm)_{[0]}=$id$_T$.
We also put $\Gam_t^\pm (a(t), b(t))=\Gam_t^\pm$.

\bgn{lemma}\label{lem: Gam-k}
The $k$ th term $(\Gam_t^\pm)_{[k]}$ is given by 
\bgn{align*}
(\Gam_t^\pm)_{[k]}= \frac 1{k!}\pi\circ(\ad_{a_k}+\ad_{b_k})\circ\h J^\pm\circ\Ad_{e^\ome}+\wtil{\Gam^\pm_k}(a_{<k}, b_{<k})
\end{align*}
where the second term $\wtil{\Gam^\pm_k}(a_{<k}, b_{<k})$ depends only on $a_1,\cdots, a_{k-1}$ and 
$b_1,\cdots, b_{k-1}$.
\end{lemma}
\bgn{proof}
Substituting the identity $\Ad_{e^{Z(t)}}=$id$ +\ad_{Z(t)}+\frac1{2!}(\ad_{Z(t)})^2+\cdots$, 
we have 
\bgn{align}
\Gam_t^\pm =&\pi\circ\Ad_{e^{Z(t)}}\circ \h J^\pm\circ\Ad_{e^\ome}\\
=&\pi\circ\(\sum_{i=0}^\infty \frac 1{i!}\ad^i_{Z(t)}\circ \h J^\pm\circ \Ad_{e^\ome}\)\\
\end{align}
Then $k$-th term is given by 
\bgn{align}
\(\Gam_t^\pm \)_{[k]}
=&\pi\circ\(\ad_{Z(t)}\circ\h J^\pm\circ\Ad_{e^\ome}\)_{[k]}+\sum_{i=2}^k\pi\circ\(\frac 1{i!}(\ad_{Z(t)}^i\circ \h J^\pm\circ\Ad_{e^\ome}\)_{[k]}\\
=&\frac 1{k!}\pi\circ(\ad_{a_k}+\ad_{b_k})\circ\h J^\pm\circ\Ad_{e^\ome}+\wtil{\Gam^\pm_k}(a_{<k}, b_{<k}),
\end{align}
where $\wtil{\Gam^\pm_k}(a_{<k}, b_{<k})$ denotes the non-linear term depending $a_1,\cdots, a_{k-1}$ and $b_1, \cdots, b_{k-1}$.
\end{proof}
\bgn{lemma}\label{lem: key b }
Let $b$ be a section of the bundle $(L^-_{\J}\cdot\ol L^+_{\J}\oplus \ol L^-_{\J}\cdot L^+_{\J})$. 
Then we have 
$$
[\,\pi(\ad_b\circ\h J^\pm\circ\Ad_{e^\ome}),\, J\, ]=0\in\text{\rm End}(T).
$$
\end{lemma}
\bgn{proof}
Applying (\ref{eq: pi0{-1}}), for $v\in T^{1,0}_J$, we obtain 
$$\h J^\pm\circ\Ad_{e^\ome} v=\Ad_{e^{\pm\sqrt{-1}\ome}}v=(\pi_0^\pm)^{-1}v\in \ol L^\pm_\J.$$
Since $\ad_b(\ol L^\pm_\J)=[b ,\ol L^\pm_\J]\subset \ol L^\mp_\J$ and $\pi (\ol L^\mp_\J)=T^{1,0}_J$,
thus we have $\pi(\ad_b\circ \h J^\pm\circ\Ad_{e^\ome} )v\in T^{1,0}_J$. 
It follows that $[\,\pi(\ad_b\circ\h J^\pm\circ\Ad_{e^\ome}),\, J\, ]=0.$
\end{proof}
The tensor space $T\otimes T^*$ defines a subbundle of CL$^2$. We denote it by $T\cdot T^*$.
An element $\gam\in T\cdot T^*$ gives the endmorphism ad$_\gam$ by ad$_\gam E=[\gam, E]$ for $E\in \TT$, which preserves the cotangent bundle $T^*$.
\bgn{lemma}\label{lem: gam}
Let $\gam$ be an element of $T\cdot T^*$. Then we have 
$$\pi\circ(\ad_{\gam}\circ \h J^\pm\circ \Ad_{e^\ome})=\ad_\gam\in \text{\rm End}(T). $$
\end{lemma}
\bgn{proof}
For a tangent vector $v\in T$, we have $\Ad_{e^\ome}v= v+[\ome, v]=v+\ad_\ome v$.
Since the map $\ad_\gam$ preserves the cotangent $T^*$, 
we have $\ad_{\gam}\circ \h J^\pm\circ \ad_{\ome}(v)\in T^*$ for all tangent $v\in T$.
Thus it follows that
$\pi( \ad_\gam\circ\h J^\pm\circ\ad_\ome)=0$, since $\pi$ is the projection to the tangent $T$.
Thus we obtain the result.
\end{proof}

\bgn{lemma}\label{lem: d Gam-k }
We assume that there is a set of sections 
$a_1,\cdots , a_k\in \CL^2$ and real sections $b_1,\cdots, b_k\in 
(L^-_{\J_0}\cdot\ol L^+_{\J_0}\oplus \ol L^-_{\J_0}\cdot L^+_{\J_0})$ which satisfies the following equations, 
\bgn{align*}
\pi_{U^{-n+3}}&\( e^{-Z(t)}\, d\, e^{Z(t)}\)_{[i]}=0, \quad 0\leq 
\forall\, i\leq k\\
&\( d e^{Z(t)}\cdot \psi_0 \)_{[i]}=0, \,\,\,\,\,\quad 0\leq 
\forall\, i\leq k\\
&[(\Gam^\pm_t)_{[i]}, \, J]=0,  \,\,\,\,\,\,\,\, \quad\quad 0\leq 
\forall\, i< k\\
\end{align*}
Then the $k$-th term $(\Gam_t^\pm)_{[k]}$ satisfies 
$$
\pi_{U^{-n+3}}[\, d,\,(\Gam_t^\pm)_{[k]}\, ]  =0,
$$
where $[\, d,\,(\Gam_t^\pm)_{[k]}\, ] $ is an operator from $U^{-n}=K_\J$ to 
$U^{-n+1}\oplus U^{-n+3}$ and $\pi_{U^{-n+3}}$ denotes the projection to the component $U^{-n+3}$.
\end{lemma}
\bgn{proof} Since we assume that the space of the obstructions to deformations of \complex structures vanishes, 
we obtain a family of section $\ch a(t)$ with $\ch a_i =a_i$ for $ i=1, \cdots k$ such that $\ch a(t)$ gives deformations of \complex structures, that is, 
$$
\pi_{U^{-n+3}}e^{-\ch a(t)} d e^{\ch a(t)} =0.
$$
The the stability theorem of \K\"ahler structures in \cite{Go1} provides deformations of \K\"ahler structures with one pure spinor, 
$(\Ad_{e^{\ch Z(t)}}\J_0, \,e^{\ch Z(t)}\psi_0)$, where $e^{\ch Z(t)}=e^{\ch a(t)}e^{\ch b(t)}$, where 
$\ch b(t)$ is a family of real sections with $\ch b_i =b_i$, for $i=1,\cdots k$.
From the correspondence between \K\"ahler structures and bihermitian structures, we  have the family of bihermitian structures 
$(J^+_t, J^-_t)$ which is given by the action of $\ch{\Gam}^\pm_t:=\Gam^\pm_t(\ch a(t), \ch b(t))$ of GL$(T)$. 
Since $J_t^\pm$ is integrable, we have 
\bgn{equation}\label{eq: ch Gam}
\pi_{U^{-n+3}}\((\ch{\Gam}^\pm_t)^{-1} \, d\,\ch{\Gam}^\pm_t\)=0.
\end{equation}
Let $\Ome$ be a $d$-closed form of type $(n,0)$ which is a local basis of $K_J=K_\J$. 
Then as in the argument of proof of the proposition \ref {prop: k th-G.complex}, we have 
$$
d\Gam^\pm_t\Ome \underset{k}\equiv \Gam^\pm_t E(t) \Ome.
$$
Since $d\Ome=0$, the degree of $E(t)$ is greater than or equal to $1$.
The condition $[(\Gam^\pm_t)_{[i]}, \, J]=0$ $( 0\leq i<k)$
implies that $(\Gam^\pm_t)_{[i]}E(t) \Ome  \in U^{-n+1}_\J$. 
Thus we have 
$$
d (\Gam^\pm_t)_{[k]}\Ome=\sum_{\stackrel {\scriptstyle i+j=k}{ 0<i,j<k}} (\Gam^\pm_t)_{[i]} E(t)_{[j]}\Ome\in U^{-n+1}_\J
$$
Hence we have 
$\pi_{U^{-n+3}}[d, (\Gam^+_t)_{[k]}]=0$.
\end{proof}
\section{Construction of deformations of bihermitian structures with $J^+_t=J$}
This section and next section are devoted to prove our main theorem
\ref{th: main theorem}.
We use the same notation as before. 
As we see in the lemma \ref{lem: Gam-k}, $(\Gam_t^\pm)_{[k]}$ depends on $a_1, \cdots, a_k$ and $b_1,\cdots, b_{k}$.  
We write $\Gam_t^\pm(a_{<k}, a_k, b_{<k}, b_k)$ for $(\Gam_t^\pm)_{[k]}$.

Let $\b$ be a holomorphic $2$-vector field on a compact K\"ahler manifold $(X, J,\ome)$, that is, $\b$ is a section of $\w^2\Theta=\w^2T^{1,0}_J$. 
For $\b$,
we shall construct a section $a(t)\in \CL^2$ and a real section $b(t)\in (L^-_{\J}\cdot\ol L^+_{\J}\oplus \ol L^-_{\J}\cdot L^+_{\J})$ such that the action of  the family of the Clifford group 
$$
e^{Z(t)} = e^{a(t)} e^{b(t)}
$$
 on $(\J, \psi)=(\J_0, \psi_0)$ gives rise to a family of \K\"ahler structures
$(\Ad_{e^{Z(t)}}\J_J, \, e^{Z(t)}\psi)$ which satisfies the following three conditions:
\bgn{align}
& \J_t:=\Ad_{e^{Z(t)}}\J_ J\,\, \text{\rm are integrable \complex structures}\label{ali: (1)}\\
&d\psi_t:=d e^{Z(t)}\psi =0 \label{ali: (2)}\\
& J^+_t = J \label{ali: (3)}, 
\end{align}
where 
$(J^+_t, J^-_t)$ denote the corresponding bihermitian structures.
It follows from 
(\ref{ali: (1)}), (\ref{ali: (2)}) that $(\J_t, \psi_t)$ are \K\"ahler structures with one pure spinor which give rise to deformations of bihermitian structures preserving $J^+_t$ from (\ref{ali: (3)}). Let 
$K_J$ be the canonical line bundle on $(X, J)$ which consists on holomorphic $n$-forms. The action
$\CL^1$ on $K_J$ provides a bundle  $\CL^1\cdot K_J$.
Then as before, the condition (\ref{ali: (1)}) is equivalent to the followings,
\bgn{equation}
{e^{-Z(t)}}\, d\,{e^{Z(t)}}\cdot K_J \subset \CL^1\cdot K_J,
\end{equation}
This implies that the ${e^{-Z(t)}}\, d\,{e^{Z(t)}}\cdot\Ome$ is written as 
$E\cdot\Ome$ for any form $\Ome$ of type $(n,0)$, where $E\in \CL^1=\TT$.
As we see in the previous section, the condition (\ref{ali: (3)}) is equivalent to 
$[\Gam^+_t(a(t), b(t) ), \, J] =0$. 
We denote by $\({e^{-Z(t)}}\, d\,{e^{Z(t)}}\)_{[i]}$ the $i$-th term of 
$\({e^{-Z(t)}}\, d\,{e^{Z(t)}}\)$ on $t$ and also write $i$-th terms of 
$d\psi_t$ and $\Gam^+_t $ by $\(d\psi_t\)_{[i]}$ and 
$\(\Gam_t^+\)_{[i]}$ respectively.
Then the three equations (\ref{ali: (1)}),(\ref{ali: (2)}) and (\ref{ali: (3)}) are reduced to the following system of equations, 
\bgn{align}
&\({e^{-Z(t)}}\, d\,{e^{Z(t)}}\)_{[i]}\cdot K_J \subset \CL^1\cdot K_J,\quad 0\leq \text{for all } i\leq k  \label{ali: (1 k)}\\
&\(d\psi_t\)_{[i]}:=(d e^{Z(t)}\psi)_{[i]} =0,\quad 0\leq \text{for all } i\leq k   \label{ali: (2 k)}\\
&[\,\(\,\Gam^+_t( a(t), b(t))\,\)_{[i]},\, J\,]=0,\quad 0\leq\text{for all } i\leq k  \label{ali: (3 k)}
\end{align}
We shall construct a solution of the system of the equations by the induction on degree $k$ of $t$. 
In the first case $k=1$, three equations are given by
\bgn{align*}
&({e^{-Z(t)}}\, d\,{e^{Z(t)}})_{[1]}\cdot K_J =[d, a_1]\cdot K_J\subset \CL^1\cdot K_J\\
&(d\psi_t)_{[1]}= d(a_1+ b_1)\cdot\psi_0 =0  \\
&[\,(\Gam^+_t(a_1,b_1))_{[1]},\,J\,]=0
\end{align*}

At first we put $\h a_1=\b+\ol\b$, where $\ol\b$ denotes the complex conjugate of $\b$. 
Since $\b$ is holomorphic, it follows that $[d, \h a_1]\cdot K_J\subset \CL\cdot K_J$.
Then from the proposition \ref{prop: GK deformations}, we have a real section $\h b_1\in (L^-_{\J}\cdot\ol L^+_{\J}\oplus \ol L^-_{\J}\cdot L^+_{\J})$ with $d (\h a_1+ \h b_1)\cdot\psi_0=0$. 
Then $\Gam_t^\pm(\h a_1, \h b_1)$ is given by 
$$
\Gam_t^\pm(\h a_1, \h b_1) =\pi\circ (\ad_{\h a_1}+\ad_{\h b_1})\circ\h J^\pm\circ\Ad_{e^\ome}.
$$
Then 
we define $\gam_1\in T\cdot T^*$ by 
\bgn{equation}\label{eq: gam 1}
ad_{\gam_1}=-\(\Gam_t^+(\h a_1, \h b_1)\)_{[1]}\in\text{\rm End}(T).
\end{equation}
It follows from the lemma \ref{lem: d Gam-k } that 
$\gam_1$ satisfies 
$[d, \gam_1]\cdot K_J \subset \CL^1\cdot K_J$, where we identify 
End$(T)$ with $T\cdot T^*$.
We define $a_1$ by 
\bgn{equation}\label{eq: a1}
a_1=\h a_1+\gam_1.
\end{equation}
Then we have 
\bgn{align}
({e^{-Z(t)}}\, d\,{e^{Z(t)}})_{[1]}\cdot K_J =&[d, a_1]\cdot K_J\\
=&\([d, \h a_1] + [d, \gam_1]\)\cdot K_J \subset \CL^1\cdot K_J
\end{align}
From the proposition \ref{prop: GK deformations}, we also have a real section 
$b_1\in (L^-_{\J}\cdot\ol L^+_{\J}\oplus \ol L^-_{\J}\cdot L^+_{\J})$ with 
$d(a_1+b_1)\cdot\psi_0=0$.
Applying the lemma \ref{lem: Gam-k} and substituting $a_1$ and $b_1$ into $\(\Gam^+_t\)_{[1]}$, we have 
\bgn{align}
(\Gam^+_t( a_1, b_1))_{[1]}=&\pi\circ(\ad_{a_1}+\ad_{b_1} )\circ\h J^+\circ\Ad_{e^\ome} \\
=&\pi\circ(\ad_{\h a_1}+\ad_{\gam_1}+\ad_{\h b_1}+\ad_{b_1-\h b_1} )\circ\h J^+\circ\Ad_{e^\ome} 
\end{align}
Applying the lemma \ref{lem: gam} to $\gam_1$ and using (\ref{eq: gam 1}), we obtain
\bgn{align}
(\Gam^+_t( a_1, b_1))_{[1]}
=&\ad_{\gam_1}+\pi\circ(\ad_{\h a_1}+\ad_{\h b_1} )\circ\h J^+\circ\Ad_{e^\ome} \\
+&\pi\circ  \ad_{b_1-\h b_1}\circ\h J^+\circ\Ad_{e^\ome} \\
=&\ad_{\gam_1}+\Gam_t^+(\h a_1, \h b_1) 
+\pi\circ \ad_{b_1-\h b_1}\circ\h J^+\circ\Ad_{e^\ome} \\
=&\pi\circ  \ad_{b_1-\h b_1}\circ\h J^+\circ\Ad_{e^\ome} 
\end{align}
 Since $b_1-\h b_1\in(L^-_{\J}\cdot\ol L^+_{\J}\oplus \ol L^-_{\J}\cdot L^+_{\J})$, it follows from the lemma \ref{lem: key b } that 
 $$
 [\ (\Gam^+_t( a_1, b_1))_{[1]},\, J\,]=0.
 $$
Hence $a_1$ and $b_1$ as above satisfies the three equations for $k=1$.
 
 We assume that there is a set of real sections $a_1,\cdots, a_{k-1}\in \CL^2$ and $ b_1,\cdots, 
 b_{k-1}\in (L^-_{\J}\cdot\ol L^+_{\J}\oplus \ol L^-_{\J}\cdot L^+_{\J})$ which satisfies the system of equations: 
 \bgn{align}
&\({e^{-Z(t)}}\, d\,{e^{Z(t)}}\)_{[i]}\cdot K_J \subset \CL^1\cdot K_J,\quad 0\leq \text{for all } i\leq k-1  \label{ali: (1 k-1)}\\
&\(d\psi_t\)_{[i]}:=(d e^{Z(t)}\psi)_{[i]} =0,\quad\quad\quad\quad 0\leq \text{for all } i\leq k-1   \label{ali: (2 k-1)}\\
&[\, \(\,\Gam^+_t( a(t), b(t))\,\)_{[i]},\,J\,]=0,\quad\quad\quad\ 0\leq\text{for all } i\leq k-1  \label{ali: (3 k-1)}
\end{align}
The $k$-th term $\({e^{- Z(t)}}\, d\,{e^{ Z(t)}}\)_{[k]}\cdot K_J$ 
is decomposed into the linear term $\frac 1{k!} [d, a_k]\cdot K_J$ and the nonlinear term $\Ob^J_k(a_{<k}, b_{<k})\cdot K_J$
which is called the term of the obstruction, 
$$
\({e^{-\h Z(t)}}\, d\,{e^{\h Z(t)}}\)_{[k]}\cdot K_J=\frac 1{k!} [d, a_k]\cdot K_J+\Ob^J_k(a_{<k}, b_{<k})\cdot K_J.
$$
We also have the decomposition of the $k$-th term $\(d e^{\h Z(t)}\psi\)_{[k]}$, 
$$
\(d e^{\h Z(t)}\psi\)_{[k]}=\frac 1{k!} d( a_k+ b_k)\psi+\Ob_k^{\psi_0}(a_{<k}, \,b_{<k})
$$
 
 From the proposition \ref{prop: GK deformations}, we have 
 the sections $\h a_k$ and 
 $\h b_k\in(L^-_{\J}\cdot\ol L^+_{\J}\oplus \ol L^-_{\J}\cdot L^+_{\J})$ which 
 satisfies 
 \bgn{align}
&\(\frac 1{k!} [d, \h a_k]+\Ob^J_k(a_{<k}, b_{<k})\)\cdot K_J\subset \CL^1\cdot K_J,\\
&\frac 1{k!} d( \h a_k+ \h b_k)\psi_0+\Ob_k^{\psi_0}(a_{<k}, \,b_{<k})=0
\end{align}
 Then from the lemma \ref{lem: Gam-k},
 $\Gam_t^+( a_{<k},\, \h a_k, \,b_{<k}, \ b_k)$ is given by 
 $$
 \Gam_t^+( a_{<k},\, \h a_k, \,b_{<k}, \ b_k)=\pi\circ(\ad_{\h a_k}+\ad_{\h b_k})\circ\h J^+\circ\Ad_{e^\ome}+\til\Gam^+_k(a_{<k}, b_{<k}).
 $$
 Then we define $\gam_k\in T\cdot T^*$ by using $\h a_k$ and $\h b_k$
 \bgn{equation}\label{eq: gam k}
 \ad_{\gam_k} =-\( \,\Gam_t^+( a_{<k},\, \h a_k, \,b_{<k}, \h b_k) \,\)_{[k]}.
 \end{equation}
 It follows from the lemma \ref{lem: d Gam-k } that 
 we have $[d ,\gam_k]\cdot K_J \subset \CL^1\cdot K_J$. 
 We define  $a_k$ by 
 \bgn{equation}\label{eq: a k}
 a_k =\h a_k +\gam_k.
 \end{equation}
 Then we have 
 \bgn{align}
 \({e^{- Z(t)}}\, d\,{e^{Z(t)}}\)_{[k]}\cdot K_J=&\(\frac 1{k!} [d, a_k]+\Ob^J_k(a_{<k}, b_{<k}) \)\cdot K_J\\
 =&\(\frac 1{k!} [d, \gam_k]+\frac 1{k!} [d, \h a_k]+\Ob^J_k(a_{<k}, b_{<k}) \)\cdot K_J\subset 
 \CL^1\cdot K_J,
 \end{align}
 where $e^{Z(t)}=e^{a(t)}\, e^{b(t)}$.
 We apply the proposition \ref{prop: GK deformations} again to obtain  a real section $b_k\in (L^-_{\J}\cdot\ol L^+_{\J}\oplus \ol L^-_{\J}\cdot L^+_{\J})$ which satisfies 
 $\( de^{Z(t)}\cdot\psi\)_{[k]}=0$.
 Then from the lemma \ref{lem: Gam-k} and (\ref{eq: a k}),
 $\(\Gam_t^+(a(t), b(t)) \)_{[k]}$ is given by 
 \bgn{align}
 k!\(\Gam_t^+(a(t), b(t)) \)_{[k]}=&
 \pi \(\ad_{a_k}+\ad_{b_k}\)\circ J^*\circ \Ad_{e^\ome} +
 k!\wtil{(\Gam_k^\pm)}(a_{<k},\, b_{<k})\\
 =&\pi\(\ad_{\h a_k}+\ad_{\gam_k}+\ad_{b_k}\)\circ J^*\circ \Ad_{e^\ome} +
 k!\wtil{(\Gam_k^\pm)}(a_{<k},\, b_{<k})\\
 \end{align}
 Applying  lemma \ref{lem: gam} to $\gam_k$ and using (\ref{eq: gam k}), we obtain
 \bgn{align}
  k!\(\Gam_t^+(a(t), b(t)) \)_{[k]}
  =&\ad_{\gam_k}+\pi\(\ad_{\h a_k}+\ad_{\h b_k}\)\circ J^*\circ \Ad_{e^\ome} \\
 +&k!\wtil{(\Gam_k^\pm)}(a_{<k},\, b_{<k})+\pi\(\ad_{b_k-\h b_k}\)\circ\h J^+\circ\Ad_{e^\ome}\\
 =&\ad_{\gam_k}+\Gam^+(a_{<k}, \h a_k, b_{<k} ,\h b_k)_{[k]}+\pi\(\ad_{b_k-\h b_k}\)\circ\h J^+\circ\Ad_{e^\ome}\\
 =&\pi\(\ad_{b_k-\h b_k}\)\circ\h J^+\circ\Ad_{e^\ome}
 \end{align}
 Since $b_k-\h b_k\in(L^-_{\J}\cdot\ol L^+_{\J}\oplus \ol L^-_{\J}\cdot L^+_{\J})$, it follows from the lemma \ref{lem: key b } that 
 $$
 [\ (\Gam^+_t( a(t), b(t)))_{[k]},\, J\,]=0.
 $$

Hence the set of sections $a_k, \, b_k$ together with $a_{<k}, b_{<k}$ satisfies 
three equations (\ref{ali: (1 k)}), (\ref{ali: (2 k)}) and (\ref{ali: (3 k)}). 
Thus from our assumption of the induction, we successively solve the equations to obtain a set of sections $a(t)$ and $b(t)$ which satisfies three equations (\ref{ali: (1 k)}), (\ref{ali: (2 k)}) and (\ref{ali: (3 k)}) for all $k$.
The solution $( a(t), b(t))$ is given in the form of a formal power series in $t$. 
Next section we shall show that both $a(t)$ and $b(t)$ are convergent series which are smooth. \\
\medskip

Our construction is well explained by the following figure,
$$
\xymatrix{
& k=1&&k=2&\\
*+[F]{\h a_1=\b+\ol\b} \ar@{->}[d] & *+[F]{a_1:=\h a_1+\gam_1} \ar@{->}[d] \ar@{=>}[r]& *+[F]{\h a_2}\ar@{->}[d] & *+[F]{a_2:=\h a_2+\gam_2}\ar@{->}[d] \ar@{=>}[r]  &\cdots\\
 *+[F]{\h b_1}\ar@{->}[d]& *+[F]{b_1}\ar@{->}[d] & *+[F]{\h b_2} \ar@{->}[d] & *+[F]{b_2} \ar@{->}[d]  &\cdots  \\
*+[F]{\gam_1=-
(\Gam^+_t(\h a_1, \h b_1))_{[1]}}\ar@{->}[ruu]& *+[F]{(\Gam^+_t(a_1, b_1))_{[1]}} & *+[F]{\gam_2=-(\Gam^+_t(a_1, \h a_2, b_1, \h b_2))_{[2]}} \ar@{->}[ruu] & *+[F]{\Gam^+_t(a(t), b(t))_{[2]}}  &\cdots \\
}
$$
\smallskip
\begin{center} Figure 1
\end{center}

\section{The convergence} 
As in the proposition \ref {prop: GK deformations}, 
if there is a set of sections $a_1,\cdots a_{k-1}$ of CL$^2$ which satisfies 
$$
\pi_{U^{-n+3}}\(e^{-a(t)}\, de^{a(t)}\)_{[i]}=0, \qquad \text{\rm for all }
i<k,
$$
and $\|a(t)\|_s\<_{k-1} K_1M(t)$,
then there is a set of real sections $b_1,\cdots, 
b_k\in (L^-_{\J}\cdot\ol L^+_{\J}\oplus \ol L^-_{\J}\cdot L^+_{\J})$    which satisfy 
the following equations:
\bgn{align}
&\pi_{U^{-n+3}}\(e^{-Z(t)}\, de^{Z(t)}\)_{[k]}=0\\
&\( d e^{Z(t)}\cdot\psi_0\)_{[i]}=0,\qquad \text{\rm for all } i\leq k\\
&\|\h a_k\|_s < K_1\lam M_k\\
&\|\h b_k\|_s <K_2M_k
\end{align}
where  $\h a_k$ is the section in the proposition \ref{prop: k th-G.complex} and
 $M(t)$ is the convergent series in (\ref{eq: M(t)}) with a constant $\lam$ and 
$K_1$ is a positive constant and a positive constant $K_2$ is determined by $\lam, K_1$.
We also have an estimate of $e^{Z(t)}=e^{a(t)}e^{b(t)}$ in \cite{Go1}, 
$$
\|Z(t)\|\< _k M(t).
$$

Then $\gam_k$ in (\ref{eq: gam k}) satisfies 
\bgn{align}
\|\gam_k\|_s <& \| \Gam^+_k (a_{<k}, \h a_k, b_{<k}, \h b_k) \|_s\\
<&2 \|\h a_k\|_s +2\|\h b_k\|_s+\| \wtil{\Gam^+_k }(a_{<k}, b_{<k})\|_s
\end{align}
Recall that $\Gam_t^+ =\pi \( \Ad_{e^{Z(t)}}\circ \h J^+\circ\Ad_{e^\ome}\)$.
Then we have an estimate of the non-linear term $\| \wtil{\Gam^+_k }(a_{<k}, b_{<k})\|_s$ 
$$\| \wtil{\Gam^+_k }(a_{<k}, b_{<k})\|_s < C\|( e^{Z(t)}-Z(t) -1)_{[k]}\|_s,$$
where $C$ denotes a constant. 
It follows from (\ref{e M(t)}) that $\|( e^{Z(t)}-Z(t) -1)_{[k]}\|_s< C(\lam)M_k$, where $C(\lam)$ satisfies 
$\lim_{\lam\to0}C(\lam)=0$. 
Thus we have 
$$
\|\gam_k\|_s< 2\|\h a_k\|_s+2\|\h b_k\|_s+ C(\lam)M_k < 2\lam K_1 M_k +
K_2M_k+C(\lam)M_k
$$
We take $\lam $ and $K_2$ sufficiently small such that $3\lam K_1M_k +K_2M_k+ C(\lam)M_k< K_1M_k$. 
Then we obtain
$$\|a_k\|_s < \|\h a_k \|_s+ \|\gam_k\|_s< K_1 M_k.$$
Thus our solution $a(t)$ satisfies that 
$\|a(t)\|_s \<_k K_1M(t)$ for all $k$.
It implies that $a(t)$ is a convergent series. 
Applying the proposition \ref{prop: GK deformations} again, we have 
$\|b(t)\|_s<_k K_2 M(t)$. 
Hence $b(t)$ is also a convergent series.
Thus it follows that $Z(t)$ is a convergent series.
\bgn{proof}{\sc of theorem \ref{th: main theorem} and theorem 0.2. }\,\,
The sections $a(t)$ and $b(t)$ which constructed in section 5 give 
deformations of bihermitian structures 
$(J^+_t, J^-_t)$.
We shall show that the family of deformations satisfies the condition in 
the theorem \ref {th: main theorem}. We already have $[\Gam^+_t, J]=0$ which implies that $J^+_t=J$. 
From the lemma \ref {lem: Gam-k} and the lemma \ref{lem: key b }, the 1st term of $J^-_t$ is given by 
\bgn{align*}
[(\Gam^-_t)_{[1]},J]=&[\pi\circ (\ad_{\h a_1} +\ad_{\gam_1} +\ad_{\h b_1})\circ \h J^-\circ\Ad_{e^\ome}), \, J]\\
=&[(\ad_{\gam_1}+\pi\circ\ad_{\h a_1}\circ\h J^-\circ\Ad_{e^\ome}),\, J]
\end{align*}
Since $\h a_1=\b+\ol\b$, we have $\pi \circ\ad_{\h a_1} |_T=0$. 
We also have ad$_{\gam_1}=-\Gam^+(\h a_1, \h b_1)$ and 
$[\ad_{\gam_1}, J] = [(\pi\circ \ad_{\gam_1}\circ J^*\circ \ad_\ome), J]$.
Thus we obtain 
$$
[(\Gam^-_t)_{[1]}, J]=2[(\pi\circ \ad_{\h a_1} \circ J^*\circ\ad_{\ome}),\, J]
$$
Then we have for a vector $v$,
\bgn{align*}
2(\pi\circ\ad_{\h a_1} \circ J^*\circ \ad_{\ome})v 
=&-2\[\b+\ol\b , \,[\ome, Jv]\]\\
=&-2\[[\b+\ol\b,\,\ome], \,Jv\]=-2(\b\cdot\ome+\ol\b\cdot\ome)Jv.
\end{align*}
Thus it follows that $\frac d{dt} J^-_t|_{t=0} =[(\Gam^-_t)_{[1]}, J]=
-2(\b\cdot\ome+\ol\b\cdot\ome)$ and the Kodaira-Spencer class of deformations $\{J^-_t\}$ is given by the class 
$-2[\b\cdot \ome]\in H^1(M.\Theta)$. 
If the Kodaira-Spencer class does not vanish, then the deformations $\{J^-_t\}$ is not trivial. 
Thus $(X, J^-_t)$ is not biholomorphic to $(X, J)$ for small $t\neq 0$. 
\end{proof}
\section{Applications}
\subsection{Bihermitian structures on del Pezzo surfaces}
A del Pezzo surface is by definition a smooth algebraic surface with ample anti-canonical line bundle.
A classification of del Pezzo surfaces are well known, they are 
$\C P^1\times \C P^1$ or $\C P^2$ or a surface $S_n$ which is the blow-up of $\C P^2$ at $n$ points $P_1,\cdots, P_n$, $(0< n \leq 8)$. 
The set of the points $\Sig:=\{P_1,\cdots, P_n\}$ must be in
{\it general position} to yield a del Pezzo surface.
The following theorem is due to Demazure, {\rm \cite{De} (see page 27)}, which shows the meaning of {\it general position},
\bgn{theorem}
The following conditions are equivalent:\\
(1) The anti-canonical line bundle of $S_n$ is ample \\
(2) No three of $\Sig$ lie on a line, no six of  $\Sig$ lie on a conic and 
no eight of $\Sig$ lie on a cubic with a double point $P_i\in \Sig$\\
(3) There is no curve $C$ on $S_n$ with $-K_{S_n}\cdot C \leq 0$.\\
(4) There is no curve $C$ with $C\cdot C=-2$ and $K_{S_n}\cdot C=0$.
\end{theorem}
\bgn{remark}
If three points lie on a line $l$, then the strict transform $\h l$ of $l$ in $S_3$ is a $(-2)$-curve with $K_{S_3}\cdot \h l=0$. 
If six points belong to a conic curve $C$, then 
the strict transform form $\h C$ of $C$ is again a $(-2)$-curve with $K_{S_6}\cdot \ol C=0$.
If eight points $P_1\cdots, P_8$ lie on a cubic curve with a double point $P_1$, then the strict transform $\h C$ of $C$ satisfies
$\h C \thicksim\pi^{-1}C -2E_1-E_2-\cdots -E_8$, where 
$E_i$ is the exceptional curve $\pi^{-1}(P_i)$. 
Then we also have $\h C^2=-2$ and $K_{S_8}\cdot \h C=0$.
\end{remark}

Let $D$ be a smooth anti-canonical divisor of ${S_n}$ which
 is given by the zero locus of a section $\b\in H^0({S_n}, K_{S_n}^{-1})$. 
Since the anti-canonical bundle $K_{S_n}^{-1}$ is regarded as the bundle of $2$-vectors $\w^2\Theta$ 
and $[\b, \b]_{S}=0\in \w^3\Theta$ on $S_n$,  every section $\b$ is  a holomorphic Poisson structure. 
On $S_n$, we have the followings,
$$
\dim H^1(S_n, \Theta)=
\bgn{cases}
&2n-8\quad ( n=5, 6,7,8) \\
&0\quad \,\,\,\qquad (n<5)
\end{cases}
$$
$$
\dim H^0(S_n, K^{-1}) =10-n
$$
and 
$$
H^{1,1}(S_n ) =1+n.
$$

Further we have 
$H^2(S_n ,\Theta)=\{0\}, \quad H^1(S_n , \w^2\Theta)\cong 
H^1(S_n , -K_{S_n})=\{0\}$.
Hence the obstruction vanishes and we have deformations of \complex structures parametrized by $H^0(S_n ,K_{S_n}^{-1})\oplus H^1(S_n , \Theta)$.

In particular,  if $n\geq 5$, we have deformations of ordinary complex structures on $S_n$
\bgn{proposition}\label{prop: del Pezzo 1}
Let $D$ be a smooth anti-canonical divisor given by the zero locus of $\b$ as above. 
Then there is a K\"aher form $\ome$ with the class 
$[\b\cdot \ome ]\neq 0\in H^1(S_n, \Theta)$.
\end{proposition}

We also have $H^2(\C P^1\times \C P^1,\Theta)=0$ and 
$H^1(\C P^1\times \C P^1, -K)=0$.

Thus we can apply our construction to every del Pezzo surface.
From the main theorem \ref {th: main theorem} together with the proposition
\ref{prop: del Pezzo 1}, we have 
\bgn{proposition}
Every del Pezzo surface admits deformations of bihermitian structures 
$(J, J^-_t, h_t)$ with $J^-_0=J$ which satisfies 
\bgn{equation}
\frac d{dt} J^-_t|_{t=0} =-2(\b\cdot\ome+\ol\b\cdot\ome),
\end{equation}
for every K\"ahler form $\ome$ and every holomorphic Poisson structure $\b$.
Further, a del Pezzo surface $S_n$ $(n\geq 5)$ admits 
distinct bihermitian structures $(J, J^-_t, h_t)$, that is, 
 the complex manifold $(X, J^-_t)$ is not biholomorphic to $(X, J)$ for small $t\neq 0 $.
\end{proposition}
Note that for small $t\neq 0$, $J^-_t\neq \pm J$.
We will give a proof of the proposition \ref{prop: del Pezzo 1} in the rest of this subsection. 

Let $N_D$ is the normal bundle to $D$ in ${S_n}$ and $i^*T_{S_n}$ the pull back of the tangent bundle $T_{S_n}$ of ${S_n}$ by the inclusion $i : D\to {S_n}$. 
Then we have the short exact sequence, 
$$
0\to T_D \to i^*T_{S_n} \to N_D\to 0
$$
and we have the long exact sequence
\bgn{align*}
0\to &H^0(D,T_D)\to 
H^0(D, i^*T_{S_n})\to 
H^0(D, N_D)\overset{\delta}\to H^1(D, T_D)\to \cdots \\
\end{align*}
Since the line bundle $N_D$ is positive, $H^1(D, N_D)=\{0\}$ and 
$\dim H^0(D, N_D)$ is equal to the intersection number $D\cdot D=9-n$ by the Riemann-Roch theorem.
Since $D$ is an elliptic curve, 
$\dim H^1(D, T_D) =\dim H^0(D, T_D) =1$. 
Hence if follows that 
\bgn{equation}\label{eq: key inequ in del Pezzo}
9-n\leq \dim H^0(D, i^*T_{S_n})\leq 10-n.
\end{equation}
Let $\I_D$ be the ideal sheaf of $D$ and $\O_D$ the structure sheaf of $D$. 
Then we have the short exact sequence
$$
0\to \I_D\to \O_{S_n}\to i_*\O_D\to 0
$$
By the tensor product, we also have 
\bgn{equation}\label{eq: exact sequence 2 in del Pezzo}
0\to \I_D\otimes T_{S_n}\to T_{S_n}\to i_*\O_D\otimes T_{S_n}\to 0
\end{equation}
Then from the projection formula we have 
$$
H^p({S_n}, i_*\O_D\otimes T_{S_n})\cong H^p({S_n}, i_*(\O_D\otimes i^*T_{S_n}))
\cong H^p(D, i^*T_{S_n}),
$$
for $p=0,1,2$.
From (\ref{eq: exact sequence 2 in del Pezzo}), we have the long exact sequence,
\bgn{align}
H^0({S_n}, T_{S_n})\to H^0(D, i^*T_{S_n})\to H^1({S_n}, \I_D\otimes T_{S_n})
\overset{j}\to H^1({S_n}, T_{S_n})\to\cdots
\end{align}
Hence we obtain
\bgn{lemma}\label{lem: key lemma in del Pezzo}
 The map $j : H^1({S_n}, \I_D\otimes T_{S_n})\to H^1({S_n}, T_{S_n})$
is not the zero map .
\end{lemma}
\bgn{proof}
We have the exact sequence,
\bgn{align}
\cdots\to H^0(D, i^*T_{S_n})\to H^1({S_n}, \I_D\otimes T_{S_n})\overset{j}\to H^1({S_n}, T_{S_n})
\end{align}
From the Serre duality with $\I_D=K_{S_n}$, we have 
$H^0(S_n, \I_D\otimes T_{S_n})\cong H^2(S_n, \Ome^1_{S_n})=\{0\}$ and 
$H^2(S_n, \I_D\otimes T_{S_n})=H^0(S_n, \Ome^1)=0$. 
From the Riemann-Roch theorem, 
$\dim H^1(S_n, \I_D\otimes T_{S_n})=n+1$.
Then it follows from (\ref{eq: key inequ in del Pezzo}) that 
$$
\dim H^0(D , i^*T_{S_n})<\dim H^1(S_n, \I_D\otimes T_{S_n})
$$
Note $10-n<n+1$ for all $n\geq 5$.
Hence the map $j$ is non zero.
\end{proof}
\bgn{remark}
Since $n\geq 5$, we have 
$H^0({S_n}, T_{S_n})=\{0\}$.
Applying the Serre duality with $K_{S_n}=\I_D$, we have $H^2(S_n, T_{S_n})\cong H^0(S_n, \I_D\otimes\Ome^1)=0$.
From the Riemann-Roch, we obtain $\dim H^1({S_n}, T_{S_n})=2n-8$.
\end{remark}
Let $\b$ be a non-zero holomorphic Poisson structure $S_n$ with the smooth divisor $D$ as the zero locus. 
Then $\b$ is regarded as a section of $\I_D\otimes \w^2\Theta$. 
Thus the section $\b\in H^0({S_n}, \I_D\otimes\w^2\Theta)$ gives an identification, 
$$\Ome^1\cong \I_D\otimes T_{S_n}.$$
Then the identification induces the isomorphism 
$$\h \b : H^1({S_n}, \Ome^1) \cong H^1({S_n}, \I_D\otimes T_{S_n}).$$
Let $j$ be the map in the lemma \ref{lem: key lemma in del Pezzo}.
Then we have the composite map $j\circ \h\b :H^1(S_n ,\Ome^1)\to H^1(S_n, \Theta)$ 
which is given by 
the class $[\b\cdot \ome] \in H^1({S_n}, T_{S_n})$ for 
$[\ome] \in H^1({S_n}, \Ome^1)$.
\bgn{proposition}
The composite map $j\circ \h \b : H^1({S_n}, \Ome^1)\to H^1({S_n}, T_{S_n})$ is not the zero map. 
\end{proposition}
\bgn{proof}
Since the map $\h\b$ is an isomorphism, $\h\b(\ome)$ is not zero.
It follows from lemma \ref{lem: key lemma in del Pezzo} that the map $j$ is non-zero.
Hence the composite map $j\circ \h \b$ is non-zero also.
\end{proof}
\bgn{proof} of lemma \ref{prop: del Pezzo 1}
The set of K\"ahler class is an open cone in 
$H^{1,1}({S_n}, \R)\cong H^2(S_n \R)$. We have the non-zero map 
$j\circ\h\b: H^{2}({S_n},\C) \cong H^1({S_n},\Ome^1) \to H^1({S_n},\Theta)$
for each $\b\in H^0({S_n}, K^{-1})$ with $\{\b=0\} =D$.
It follows that the kernel $j\circ\h\b$ is a closed subspace and 
 the intersection $\ker(j\circ\h\b)\cap H^{2}({S_n},\R)$ is closed in 
 $H^{2}({S_n}, \R)$ whose dimension is strictly less than $\dim H^{2}({S_n}, \R)$.
 Thus the complement in the K\"ahler cone 
 $$\big\{ [\ome] :\text{\rm K\"ahler class} \, |\, j\circ\h\b([\ome])\neq0\big\}$$
 is not empty. 
 Thus there is a K\"ahler form $\ome$ such that the class $[\b\cdot\ome]\in H^1(S_n ,\Theta)$ does not vanish for $n\geq 5$.
\end{proof}
We also remark that our proof of the lemma \ref {prop: del Pezzo 1} still works for degenerate del Pezzo surfaces. 
  \subsection{Vanishing theorems on surfaces}
 Let $M$ be a compact complex surface with canonical line bundle $K_M$.
  We shall give some vanishing theorems of the cohomology groups
  $H^1(M, -K_M)$ and $H^2(M, \Theta)$ on a compact smooth complex surface $M$, which are the obstruction spaces to deformations of \complex structures starting from the ordinary one $(X, \J_J)$.
The following is  practical to show the vanishing of $H^1(M, -K_M)$. 
 \bgn{proposition}\label{vanishing H1(-K)}
  Let $M$ be a compact complex surface with $H^1(M, {\cal O}_M)=0$.
  If $-K_M=m[D]$ for a irreducible, smooth curve $D$ with positive self-intersection number $D\cdot D>0$  and a positive integer $m$, then 
  $H^1(M, K^n_M)=0$ for all integer $n$.
  \end{proposition}
  The proposition is often used in the complex geometry. For completeness, we give a proof.
  \bgn{proof}
  Let $I_D$ be the ideal sheaf of the curve $D$. Then we have the short exact sequence, $0\to I_D\to {\cal O}_M\to j_*{\cal O}_D\to 0$, where $j: D\to X$.
 Then we have the exact sequence, 
 $$
 H^0(M, {\cal O}_M)\to H^0(M, j_*{\cal O}_D)\overset{\del}{\to}H^1(M, I_D)\to 
 H^1(M, {\cal O}_M) 
 $$
 It follows that the coboundary map $\del$ is a $0$-map. Thus  from $H^1(M, {\cal O}_M)=0$, we have $H^1(M,I_D)=H^1(M, -[D])=0$.
 We use the induction on $k$. We assume that $H^1(M, I_D^k)=H^1(M, -k[D])=0$ for a positive integer $k$.
 The short exact sequence $0\to I^{k+1}_D\to I^k_D\to j_*{\cal O}_D\otimes I^k_D\to 0$ induces the exact sequence, 
 $$
 H^0(M, j_*{\cal O}_D\otimes I_D^k)\to H^1(M, I_D^{k+1})\to H^1(M, I_D^k).
 $$
 By the projection formula, we have $H^0(M, j_*{\cal O}_D\otimes I_D^k)=H^0(D, -k[D]|_D)$. Since $D\cdot D>0$, it follows that the line bundle $ -k[D]|_D$ is negative and then
 $H^0(D,  -k[D]|_D)=H^0(M, I_D^k)=0$. It implies that $H^1(M, I_D^{k+1})=H^1(M, -(k+1)[D])=0$. 
 Thus by the induction, we have $H^1(M, -nD)=0$ for all positive integer $n$. 
 Applying the Serre duality, we have $H^1(M, -nD)\cong H^1(M, (n-m)D)=0$. 
 Thus $H^1(M, nD)=0$ for all integer $n$. 
 Then the result follows since $H^1(M, K^n)=H^1(M, -(nm)D)=0$.
  \end{proof}
 The author also refer to the standard vanishing theorem.
 If $D=\sum_i a_i D_i$ is a $\Q$-divisor on $M$, where $D_i$ is a prime divisor and $a_i \in \Q$. Let $\lceil a_i\rceil$ be the round-up of $a_i$ and $\lfloor a_i\rfloor$ the round-down of $a_i$. 
 Then the fractional part $\{a_i\}$ is $a_i -\lfloor a_i\rfloor$. 
 Then the round-up and the round-down of $D$ is defined by 
 $$
 \lceil D\rceil =\sum_i \lceil a_i\rceil D_i, \quad 
 \lfloor D\rfloor=\sum_i \lfloor a_i \rfloor D_i
 $$
and  $\{ D\}=\sum_i\{a_i\}D_i$ is the fractional part of $D$. 
A divisor $D$ is {\bf nef} if one has $D\cdot C\geq 0$ for any curve $C$. 
A divisor $D$ is {\bf nef} and {\bf big} if in addition, one has $D^2>0$.
We shall use the following vanishing theorem. 
The two dimensional case is due to Miyaoka and the higher dimensional cases
are due to Kawamata and Viehweg
\bgn{theorem} \label{th: vanishing K&V}
Let $M$ be a smooth projective surface and $D$ a $\Q$-divisor on $M$ such that \\
(1) supp$\{D\}$ is a divisor with normal crossings, \\
(2) $D$ is nef and big.\\
Then $H^i(M, K_M+\lceil D\rceil) =0$ for all $i>0$.
\end{theorem}
If $-K_M=mD$ is nef and big divisor where $D$ is smooth for $m>0$.
Then applying the theorem, we have 
$$
H^i(M, -K_M) \cong H^i(M, K_M-2K_M)=0,
$$
for all $i>0$.
\\
Next we consider the vanishing of the cohomology group $H^2(M, \Theta)$. 
Applying the Serre duality theorem, 
we have 
$$
H^2(M, \Theta)\cong H^0(M, \Ome ^1\otimes K_M)
$$
If $-K_M$ is an effective divisor $[D]$, then $K_M$ is given by the ideal sheaf 
$I_D$ of $D$. 
The short exact sequence: $0\to \Ome^1\otimes I_D \to \Ome^1\to 
\Ome^1\otimes{\cal O}_D\to 0$ gives us the injective map, 
$$
0\to H^0(M, \Ome^1\otimes K_M)\to H^0(M, \Ome^1).
$$ 
Hence we have 
\bgn{proposition}\label{vanishing H2Theta}
if $M$ is a smooth surface with effective anti-canonical divisor satisfying
$H^0(M, \Ome^1)=0$, then we have the vanishing $H^2(M, \Theta)=0$.
\end{proposition}

 \subsection{Non-vanishing theorem}
\bgn{proposition}\label{non-vanishing}
 Let $M$ be a K\"ahler surface with a K\"ahler form $\ome$ and 
 a non-zero Poisson structure $\b \in H^0(M, \w^2\Theta)$. 
 Let $D$ be the divisor defined by the section $\b$. 
 If there is a curve $C$ of $M$ with $C\,\cap$ supp $D=\emptyset$, 
 then the class $[\b\cdot\ome]\in H^1(M, \Theta)$ does not vanish.
 \end{proposition}
 \bgn{proof}
 Since $\b$ is not zero on the complement $M\bsh D$, there is a holomorphic symplectic form $\h \b$ on the complement.
 The symplectic form $\h \b$ gives the isomorphism 
$\Theta \cong \Ome^1$ on $M\bsh D$ which induces the isomorphism 
between cohomology groups
$H^1(M\bsh D, \Theta)\cong H^1(M\bsh D, \Ome^1)$.
Then the restricted class $[\b\cdot\ome] |_{M\bsh D}$ corresponds to the 
K\"ahler class $[\ome]|_{M\bsh D}\in H^1(M\bsh D, \Ome^1)
\cong H^{1,1}(M\bsh D)$ under the isomorphism. 
Since there is the curve $C$ on the complement $M\bsh D$ and $\ome$ is a K\"ahler form, the class $[\ome|_C]\in H^{1,1}(C)$ does not vanish.
Then it follows that the class $[\ome]|_{M\bsh D}\in H^1(M\bsh D, \Ome^1)$ does not vanish. 
It implies that $[\b\cdot\ome]|_{M\bsh D}$ does not vanish also.
Thus we have that the class $[\b\cdot \ome]\in H^1(M, \Theta)$ does not vanish.
 \end{proof}
\subsection{Deformations of bihermitian structures on the Hirtzebruch surfaces $F_2$ and $F_3$} 
Let $F_2$ be the projective space bundle of 
$T^*\C P^1\oplus {\cal O}_{\C P^1}$, 
$$
F_2=\Bbb P(T^*\C P^1\oplus {\cal O}_{\C P^1}).
$$

We denote by $E^+$ and $E^-$ the sections of $F_2$ with positive and negative self-intersection numbers respectively. 
An anti-canonical divisor of $F_2$ is given by $2E^+$, 
while the section $E^-$ with  
 $E^-\cdot E^-=-2$ is the curve which satisfies $E^+\cap E^-=\emptyset$.
 Thus we have the non-vanishing class $[\b\cdot \ome]\in H^1(F_2, \Theta)$, where $\b$ is a section of $-K$ with the divisor $2E^+$.
  \
(Note that the canonical holomorphic symplectic form $\h \b$ on the cotangent bundle $T^*\C P^1$ which induces the holomorphic Poisson structure $\b$. 
The structure $\b$ can be extended to $F_2$ which gives the anti-canonical divisor $2[E^+]$.)
\bgn{proposition}\label{prop: non-vanishing on F2}
The class $[\b\cdot\ome] \in H^1(F_2, \Theta)$ does not vanish for every K\"ahler form $\ome $ on $F_2$.
\end{proposition}
\bgn{proof}
The result follows from the proposition \ref{non-vanishing}.
\end{proof}

On the surface $F_2$, 
 the anti-canonical line bundle of $F_2$ is 
$2E^+$ and $H^1(F_2,{\cal O}_{F_2})=0$.
Hence from the proposition \ref{vanishing H1(-K)}, we  have the vanishing $H^i(F_2, -K_X)=\{0\}$ for all 
$i>0$. 
Since the surface $F_2$ is simply connected, it follows from the proposition 
\ref{vanishing H2Theta} that
 $H^2(F_2, \Theta)=0$. 
Hence the obstruction vanishes and we can apply our main theorem . 
It is known that every non-trivial small deformation of $F_2$ is $\C P^1\times \C P^1$. Thus we have
\bgn{proposition}
Let $(X, J)$ be the Hirtzebruch surface $F_2$ as above.
Then there is a family of deformations of bihermitian structures $(J^+_t, J^-_t, h_t)$
with  $J^+_t=J^-_0=J$ such that 
$(X, J^-_t)$ is $\C P^1\times \C P^1$ for small $t\neq 0$.
\end{proposition}

Let $F_e$ be the projective space bundle  
 $\Bbb P({\cal O}\oplus {\cal O}(-e))$ over $\C P^1$
 with $e>0$. 
 There is a section $b$ with $b^2=-e$, which is unique if $e>0$. 
 Let $f$ be a fibre of $F_e$. 
 Then $-K$ is given by $2b+(e+2)f$, which is an effective divisor.
Thus from the proposition \ref{vanishing H2Theta}, 
we have $H^2(F_e, \Theta)=\{0\}$. 
$P_{-1}(F_e)=\dim H^0(F_e, K^{-1})$ is listed in the table 7.1.1 of \cite{Sa1}, 
$$
P^{-1}(F_e) =
\bgn{cases}
&9 \qquad\,\,\,\, e=0,1 \\
&9 \qquad\,\,\,\, e=2\\
&e+6\quad e\geq 3
\end{cases}
$$
Since $K$ is given by the ideal sheaf $I_D$ for the effective divisor $D=2b+(e+2)f$,
It follows from the Serre duality that $H^2(F_e, K^{-1}) =H^0(F_e, I_D^2)=\{0\}$.
Thus applying the Riemann-Roch theorem, we obtain
$$
\dim H^1(F_e, K^{-1}) =e -3,
$$
for $e\geq 3$. 
In the case $e=3$, we have $H^1(F_3, K^{-1})=H^2(F_3, \Theta)=\{0\}$. Thus from the theorem \ref{th: main theorem}, we have 
\bgn{proposition}
The Hirtzebruch surface $F_3$ admits deformations of bihermitian structures 
$(J, J^-_t, h_t)$ with $J^-_t\neq \pm J$ for small $t\neq 0$.
\end{proposition}

We can generalized our discussion of $F_2$ to
 the projective space bundle of $T^*M\oplus{\cal O}_M$
 over a compact K\"ahler manifold $M$. 
 Then we also have the Poisson structure $\b$  and as in the proposition \ref{non-vanishing}, it is shown that the class $[\b\cdot \ome]$
 does not vanish. 
 Thus we have the deformations of bihermitian structures from 
 the stability theorem \cite{Go2}.
 For the ones as in the theorem \ref{th: main theorem}, we need to show the vanishing of the obstruction. 
 Note that the obstruction space does not vanish in general. 
 
 If $M$ is a Riemannian surface $\Sig_g$ of genus $g\geq 1$, then 
 the projective space bundle is called a ruled surface of degree $g$. 
 it is known that small deformations of any ruled surface of degree $g \geq 1$ 
 remain to be ruled surfaces of the same degree. 
 Applying the stability theorem, we have
\bgn{proposition}
Let $(X, J)$ be a ruled surface $\Bbb P(T^*\Sig_g\oplus {\cal O}_{\Sig_g})$
with degree $g\geq 1$. 
Then there is a family of non-trivial bihermitian structures 
$(J^+_t, J^-_t, h_t)$  such that $J\neq \pm J^\pm_t$ and 
$(X, J^\pm_t)$ is a ruled surface for small $t$.
\end{proposition}
\subsection{Bihermitian structures on degenerate del Pezzo surfaces}
We shall consider the blow-up of $\C P^2$ at $r$ points which are not in 
general position. 
We follow the construction as in 
\cite {De}, (see page 36). 
 We have a finite set $\Sig=\{x_1,\cdots, x_r\}$ and $X(\Sig)$ obtained by successive blowing up at $\Sig$,
$$
X(\Sig)\to X(\Sig_{r-1})\to \cdots \to X(\Sig_1)\to \C P^2,
$$
At first $X(\Sig_1)$ is the blow-up of $\C P^2$ at a point $x_1\in \C P^2$
 and  we have $\Sig_i =\{x_1, \cdots, x_{i}\}$ and 
$X(\Sig_{i+1})$ is the blow-up of $X(\Sig_{i})$ at $x_{i+1}\in X(\Sig_i)$.
Let $E_i$ be the divisor given by  the inverse image of $x_i\in X(\Sig_{i-1})$.
If $\Gam$ is an effective divisor on $\C P^2$, one notes that mult$(x_i, \Gam)$ 
the multiplicity of $x_i$ on the proper transform of $\Gam$ in $X(\Sig_{i-1})$, 
and one says that $\Gam$ passes through $x_i$ if mult$(x_i, \Gam)>0$. 
 Define $\h E_1, \cdots, \h E_r$ by recurrence as follows, 
 On $X(\Sig_1)$, one put $\h E_1=E_1$ ; 
 on $X(\Sig_2)$, $\h E_1$ is a proper transform of the previous $E_1$ and 
 one also put $\h E_2=E_2$; on $X(\Sig_3)$, $\h E_1$ and $\h E_2$ are the proper transform of previous $\h E_1$ and $\h E_2$  respectively and $\h E_3=E_3$. 
 Then $\h E_1, \cdots, \h E_r$ are irreducible components of $E_1+\cdots+ E_r$. 
 
We assume that the following condition on $\Sig$, \medskip\\
(*)  For each $i=1,\cdots, r$, a point $x_i\in X(\Sig_{i-1})$ does not belong to a irreducible curve $\h E_j$ with self-intersection number $-2$ for $1\leq j\leq i-1$.\\ \\
If a point $x_i\in X(\Sig_{i-1})$ belongs to a irreducible curve 
$\h E_j$ with self-intersection number $-2$, then the proper
transform of $\h E_j$ becomes a curve with self-intersection number $-3$.
If there is a rational curve with self-intersection number $-3$ or less, the anti-canonical divisor of $X(\Sig)$ is not nef.
 
\bgn{definition} 
 A set of points $\Sig $ is in {\it almost general position} if $\Sig$ satisfies the following:\\
 (1) $\Sig$ satisfies the condition (*)   \\
 (2) No line passes through $4$ points of $\Sig$ \\
 (3) No conic passes through $7$ points of $\Sig$ \\
 \end{definition}
 We call $X(\Sig)$ {\it a degenerate del Pezzo surface} if $\Sig$ is  in almost general position.
Note that if $\Sig$ is in general position, $\Sig$ is in almost  general position. 
In \cite{De}, the following theorem was shown,\\
\bgn{theorem}{\rm \cite{De} }The following conditions are equivalent: \\
(1) $\Sig$ is in almost general position \\
(2) The anti-canonical class of $X(\Sig)$ contains a smooth and irreducible curve $D$. \\
(3) There is a smooth curve of $\C P^2$ passing all points of $\Sig$.\\
(4) $H^1(X(\Sig), K_{X(\Sig)}^n)=\{0\}$ for all integer $n$\\
(5) $-K_{X_{\Sig}}\cdot  C\geq 0$ for all effective curve $C$ on $X(\Sig)$ and in adition,
 if $-K_{X(\Sig)}\cdot C =0$, then $C\cdot C=-2$.
\end{theorem}
Then from (2) there is a smooth anti-canonical divisor on a degenerate del Pezzo surface and we have $H^1(X(\Sig),{\cal O}_X)=0$.
Hence from the proposition \ref{vanishing H1(-K)}, we have the vanishing $H^i(X(\Sig), -K_X)=0$, for all $i>0$. 
  A degenerate del Pezzo surface $X(\Sig)$ satisfies $H^0(X(\Sig),\Ome^1)=0$. 
 Then it follows from the proposition \ref{vanishing H2Theta} that $H^2(X(\Sig), \Theta)=0$.
 
Let $X(\Sig)$ be a degenerate del Pezzo surface which is not a del Pezzo surface, that is,  the anti-canonical class of $X(\Sig)$ is not ample. 
Then from (5), there is a $(-2)$-curve $C$ with $K_{X(\Sig)}.\cdot C=0$. 
Then it follows that $C$ is a $\C P^1$.
Thus we contract $(-2)$-curves on a degenerate del Pezzo to obtain a complex surface with rational double points, which is called the Gorenstein log del Pezzo surface. 
Let $\b$ be a section of $-K_X$ with the smooth divisor $D$ as the zero set. 
We denote by $J$ the complex structure of the del Pezzo surface $X(\Sig)$.
From the theorem \ref{th: main theorem}, we have 
\bgn{theorem}
A degenerate del Pezzo surface admits deformations of distinct bihermitian structures $(J, J^-_t, h_t)$ with $J_0^-=J$ and $J^-_t\neq \pm J$ for 
small $t\neq 0$, that is, 
 $\frac d{dt} J^-_t|_{t=0}=
-2(\b\cdot \ome+\ol\b\cdot\ome)
$,
and 
the complex structure $J^-_t$ is not equivalent to 
$J$ of $X(\Sig)$ under diffeomorphisms for small $t\neq 0$, 
where $\ome$ is a K\"ahler form.
\end{theorem}
 \bgn{proof}
 If $X(\Sig)$ is a del Pezzo surface, we already have the result.
 If $X(\Sig)$ is not a del Pezzo but a degenerate del Pezzo, we still have $H^2(X(\Sig), \Theta)=
 H^1(X(\Sig), K^{-1})=\{0\}$. Thus we have deformations of bihermitian structures as in 
 the theorem \ref{th: main theorem}. 
 It is sufficient to show that the class $[\b\cdot\ome]$ does not vanish.
 Since $K\cdot C=0$, the line bundle $K|_C\to C\cong \C P^1$ is trivial. 
 If there is a point $P\in D\cap C$, then $\b(P)=0$ and it follows that $\b|_C\equiv 0$. Since $D$ is smooth, we have $D=C$. 
 However $D\cdot D=9-r$ and $D\cdot C= -K\cdot C=0$. 
 Thus $D\cap C=\emptyset$. 
 Then applying the proposition \ref{non-vanishing}, 
 we obtain $[\b\cdot \ome]\neq 0\in H^1(X(\Sig), \Theta)$.
 \end{proof} 
 

\section{Appendix I (The Kuranishi family of \complex structures)}
We shall discuss an analog of the Kuranishi family of 
deformations of \complex structures.
The deformation theory of \complex structures was already obtained in \cite{Gu1} by using the implicit function theorem. 
For the completeness of this paper, we will give the different construction of deformations of \complex structures by using the power series.
Our method  explicitly shows that 
the deformations family depends holomorphically on the parameter $t$ and we can also have an estimate of the convergent series as in section 1. \par
Let $(X, \J)$ be a compact \complex manifold and 
$\oL:=\ol L_\J$ the Lie algebroid bundle as before which gives the decomposition, 
${(\TT)}^\C=L\oplus \oL$.
Note that the obstruction space 
$H^3(\w^\bullet\ol L_\J)$ does not necessary vanish. 
Even in the case we obtain the family of deformations which is parametrized by 
an analytic set. We fix a metric on $X$ and consider the adjoint $d^*_L$, where $d_L$ is the derivative of the complex, 
$$
\cdots\overset{d_L}\to \w^k \ol L_\J\overset{d_L}\to \w^{k+1}\ol L_\J
\overset{d_L}\to\cdots.
$$
We also denote by $G_L$ the Green operator of the Laplacian $\trian_L:= d_Ld_L^*+ d_L^* d_L$. 

Let $\{\eta_i\}_{i=1}^m$ be a basis of the Harmonic forms $\Bbb H^2(\oL)
\cong H^2(\oL)$.
As in (\ref{eq: e(t)})
we also have the convergent series 
$\e(t)$ which is a unique solution of  
\bgn{equation}\label{eq: e(t) Kuranishi}
\e(t)=\e_1(t) -\frac12 d_L^* G_L [\e(t),\e(t)]_S,
\end{equation}
where $\e_1(t)=\sum_{i=1}^m \eta_i t_i$ and 
$t=(t_1,, \cdots, t_m )\in \C^m$.
Note that $\e(t)$ is not a section with one variable but one with several variables $t=(t_1,, \cdots, t_m )$. 
The convergent series $\e(t)$ is determined by the first term $\e_1(t)$. 
The harmonic component of $[\e(t), \e(t)]_S$ is denoted by 
$H([\e(t), \e(t)]_S)\in \H^3(\w^\bullet\oL)$. 
We define an analytic set $A$ by 
$$
A= \{ \, t\in \C^m\,\,\, \big|\, \,|t|< \a, \, H([\e(t), \e(t)]_S)=0\, \}
$$
where $\a$ is a sufficiently small constant.
 \bgn{proposition}
We have a family of \complex structures $\{\J_t\}$ which is parametrised by 
the analytic set $A$.
\end{proposition}
Our proof is almost same as in the one of complex deformations and we use the similar notation as in \cite{Ko}.
\bgn{proof}
It suffices to show that for a fixed $\e_1(t)$, the $\e(t)$ in
(\ref{eq: e(t) Kuranishi}) satisfies the Maurer-Cartan equation if and only if 
$H([\e(t), \e(t)]_S)=0$.
If $\e(t)$ is a solution of the Maurer-Cartan equation, 
$$
d_L\e(t)+\frac12[\e(t), \e(t)]_S=0.
$$
Then it follows that the harmonic part $H([\e(t), \e(t)]_S)$ vanishes.
Conversely, we assume that $H([\e(t), \e(t)]_S)=0$.
Let $\Psi=d_L\e(t)+\frac12[\e(t),\e(t)]_S\in \w^3\oL$. 
It follows from (\ref{eq: e(t) Kuranishi}) that that 
$d_L\e(t)=-\frac12d_Ld_L^*G_L[\e(t), \e(t)]_S$.
Then applying the Hodge decomposition to $[\e(t), \e(t)]_S$, we have
\bgn{align}
2\Psi=&-d_Ld_L^*G_L[\e(t), \e(t)]_S+[\e(t), \e(t)]_S\\
=&H([\e(t), \e(t)]_S)+d_L^*d_LG_L[\e(t), \e(t)]_S\\
=&d_L^*d_LG_L[\e(t), \e(t)]_S
\end{align}
By using the proposition \ref{prop: Schouten relations} and substituting 
$d_L\e(t)=\Psi-\frac12[\e(t),\e(t)]_S$, 
we have 
\bgn{align}
\Psi=&d_L^*G_L[d_L\e(t), \e(t)]_S\\
=&d_L G_L [\Psi, \e(t)\,]_S -d_L G_L \frac12\[ [\e(t),\e(t)]_S, \,\e(t)\]_S\\
=&d_L G_L [\Psi, \e(t)\,]_S.
\end{align}
We use the Sobolev norm $\|\, \|_s$ and the elliptic estimate, 
\bgn{align}
\|\Psi\|_s<& C_1\|\,[\Psi, \e(t)\,]_S\,\|_{s-1}\\
&<C_2 \|\Psi\|_s\,\|\e(t)\|_s,
\end{align}
where $C_1, C_2$ are positive constants. 
Thus for small $t$ such that $C_2\|\e(t)\|_s<1$, it follows that 
$\Psi=0$. 
Hence $\e(t)$ satisfies the Maurer-Cartan equation.
\end{proof}

\section{Appendix II}
We will give a short explanation of the Schouten bracket and the proposition \ref{prop: Schouten relations}.
Our definition of the Schouten bracket is called the Dervied bracket construction \cite{KoSc}.
Let $(X, \J)$ be a \complex manifold with th decomposition 
$(\TT)^\C=L_\J\oplus\ol L_\J$. 
We denote by $\w^\bullet\ol L_\J$ the skew-symmetric forms of $\ol L_\J$, which acts on differential forms $\w^\bullet T^*$ by the spin representation. Let $K_\J$ be the canonical line bundle which is given by 
$K_\J =\{\, \phi \in \w^\bullet T^*\, |\, L_\J \cdot \phi=0\, \}$.
Then the space of differential forms is decomposed into irreducible representations: 
$\w^\bullet T^* =\oplus_{i=0}^{2n} U^{-n+p}$, where each component $U^{-n+p}$ is given by $\w^p\ol L_\J \cdot K_\J$.
  For a section $\e\in\w^p\ol{L}_\J$, we denote by $|\e|:=p$ the degree of $\e$.
The exterior derivative  $d$ is decomposed into 
$d=\pa+\ol\pa$, where $\pa: U^{-n+p}\to U^{-n+p-1}$ and the complex conjugate 
$\ol\pa : U^{-n+p}\to U^{-n+p+1}$. 
We consider $\e\in \w^\bullet\ol L_\J$ is an operator from $K_\J $ to 
$U^{-n+|\e|}$ by the spin representation of $\w^\bullet \ol L_J$ on 
$\w^\bullet T^*$.
For $\e_1,\e_2\in \w^\bullet\ol{L}_\J$ , we define a graded bracket $[\, ,\,]_G$ by
$[\e_1, \e_2]_G=\e_1\e_2-(-1)^{|\e_1|\, |\e_2|}\e_2\e_1$.
Let $A$ be a differential operator acting on 
$\w^\bullet T^*$. If 
$A :U^{-n+i}\to U^{-n+i+a}$, for all $i$, $A$ is an operator of
degree $a=|A|$. For operators
$A,B$ of degree $|A|$ and $|B|$, 
we also have the graded bracket: 
$$[A,B]_G:=AB-(-1)^{|A|\,|B|}BA.$$

The exterior derivative $d$ admits the decomposition $d=\pa+\ol\pa$, where
$\pa$ and $\ol\pa$ are operators of degree $1$ and $-1$ respectively.
Then $d$ is an operator of odd degree and the graded commutator with
$\e\in\w\ol{L}_\J$ is given by 
$$D\e:= [d, \e]_G=d\e-(-1)^{|\e|}\e d.$$ Then we define Schouten bracket $[\e_1,,\e_2,]_S\in\w^{|\e_1|+|\e_2|-1}\ol{L}_\J$ by 
\bgn{equation}\label{eq: Schouten bracket}
[\e_1, \e_2]_S:=[\, D\e_1, \, \e_2\,]_G=[\, [d, \e_1]_G,\, \e_2\,]_G
=[\, [\pa, \e_1]_G,\, \e_2\,]_G,
\end{equation}
where $[\, [\ol\pa, \e_1]_G,\, \e_2\,]_G=0$ and 
$[\e_1,,\e_2,]_S\in \w^{|\e_1|+|\e_2|-1}\ol L_\J$.
Let $d_L$ be the derivative of the Lie algebroid $\ol L_\J$.
Then we have 
\bgn{equation}\label{eq: dL }
d_L\e =[\ol\pa, \e]_G \in \w^{p+1}\ol{L}_\J
\end{equation}
(Refer to \cite{Go0}.)
In fact, since  we have  $[\ol\pa, \e]_Gf\phi=f[\ol\pa, \e]_G+[\ol\pa f, \e]_G\phi=f[\ol\pa, \e]_G$, 
the operator  $[\ol\pa, \e]_G$ is regarded as an element of Hom$(K_\J, U^{-n+|\e|+1})$
 and since $[[\ol\pa, \e]_G,\,\e_1]_G\phi=0$ for 
$\phi\in K_\J$ and $\e_1\in \w^\bullet \ol L_\J$, the commutator
$[\ol\pa, \e]_G $ is also an element of  $\w^{|\e|+1}\ol{L}_\J$ under the isomorphism
 $\w^{p}\ol{L}_\J\cong $Hom$(K_\J, U^{-n+p})$, which is given by the spin representation.
Then we obtain an isomorphism between two complexes: 
$$
(\w^\bullet\ol{L}_\J, d_L)\cong ( U^{-n+\bullet}\otimes K^{-1}_\J, 
[\ol\pa, \,\,\,]_G)
$$
In fact we have 
\bgn{align}
\big[\,\ol\pa,\, [\ol\pa, \e]_G\,\big]_G=&\big[\,\ol\pa, \, (\ol\pa\e-(-1)^{|\e|}\e\ol\pa)\,\big]_G\\
=&\ol\pa\ol\pa\e-(-1)^{|\e|}\ol\pa\e\ol\pa-(-1)^{|\e|+1}
\ol\pa\e\ol\pa-\e\ol\pa\ol\pa\\
=&0
\end{align}
From now we identify $d_L\e$ with $[\ol\pa, \e]_G$.

We have the following relations of the graded bracket.
\bgn{lemma}\label{lem: graded bracket}
\bgn{align}
&[A, B]_G=-(-1)^{|A|\, |B|}[B,A]_G, 
\end{align}
 the  Jacobi identity of the graded bracket holds
\bgn{align}
&\big[\, [A,\,B]_G,\, C\big]_G (-1)^{|A||C|}
+\big[\,[ B, \, C]_G,\,A\big]_G(-1)^{|B||A|} +
\big[\,[C,\, A]_G,\, B\big]_G(-1)^{|C||B|}=0
\end{align}
\end{lemma}
\bgn{proof}
These follows from a direct calculations.
\end{proof}
We also have the following three relations of the Schouten bracket, 
\bgn{lemma}\label{lem: S-1}
$$[\e_1,\e_2]_S=(-1)^{|\e_1||\e_2|}[\e_2,\e_1]_S$$
\end{lemma}
\bgn{lemma}\label{lem: S-2}
$$d_L[\e_1, \e_2]_S=[d_L\e_1,\e_2]_S+(-1)^{|\e_1|}\,[\e_1, d_L\e_2]_S
$$
\end{lemma}
\bgn{lemma}\label{lem: S-3}
$$
\big[\,\,[\e_1,\e_2]_S, \,\e_3\,\big]_S(-1)^{|\e_1||\e_3|}+
\big[\,\,[\e_2,\e_3]_S, \,\e_1\,\big]_S(-1)^{|\e_2||\e_1|}+
\big[\,\,[\e_3,\e_1]_S, \,\e_2\,\big]_S(-1)^{|\e_3||\e_2|}=0
$$
\end{lemma}
We shall show that every lemma follows from (\ref{eq: dL }) and lemma \ref{lem: graded bracket}.
 \par
{\it Proof of lemma \ref{lem: S-1}}
for $\e_1, \e_2\in\w^\bullet\ol{L}_\J$, we have 
\bgn{equation}\label{eq: app D}
D[\e_1, \e_2]_G=[D\e_1, \e_2]_G+(-1)^{|\e_1|}[\e_1, D\e_2]_G=0
\end{equation}
Since $[\e_1,\e_2]_G=0$, we have 
$[D\e_1, \e_2]_G+(-1)^{|\e_1|}[\e_1, D\e_2]_G=0.$
Since $[\e_1, D\e_2]_G=-(-1)^{|\e_1|\,(|\e_2|+1)}[D\e_2,\e_1]_G$, 
we obtain 
$$
[D\e_1, \e_2]_G=(-1)^{|\e_1|\,|\e_2|}[D\e_2, \e_1]_G
$$
It implies that $[\e_1, \e_2]_S=(-1)^{|\e_1|\,|\e_2|}[\e_s, \e_1]_S$.
\qed\par
{\it Proof of lemma \ref{lem: S-2}}\,\, From (\ref{eq: Schouten bracket}) we have 
\bgn{align}
d_L[\e_1, \e_2]_S=&\[\ol\pa, \,\,[\e_1, \e_2]_S\]_G
=\[\ol\pa, \,\,\[D\e_1, \e_2]_G\]_G
\end{align}
Applying the lemma \ref{lem: graded bracket}, we have  
\bgn{align}
(-1)^{|\e_2|}d_L[\e_1, \e_2]_S=&
\[[D\e_1, \e_2]_G,\,\ol\pa\]_G(-1)(-1)^{|\e_1|+|\e_2|-1}(-1)^{|\e_2|}\\
=&\[[D\e_1, \e_2]_G,\,\ol\pa\]_G(-1)(-1)^{|\e_1|-1}\\
=&\[[\e_2, \ol\pa]_G,\,D\e_1\]_G(-1)^{|\e_2|\,(|\e_1|-1)}\\
+&\[[\ol\pa, D\e_1]_G,\,\e_2\]_G(-1)^{|\e_2|}
\end{align}
 From the lemma \ref{lem: graded bracket} , we also have 
\bgn{align}
[\ol\pa, D\e_1]_G=&[D\e_1, \ol\pa]_G(-1)(-1)^{|\e_1|-1}\\
=&\[[\pa,\e_1]_G,\,\ol\pa\]_G(-1)(-1)^{|\e_1|-1} \\
=&\[[\e_1,\ol\pa]_G,\,\pa\]_G(-1)^{|\e_1|-1}(-1)^{|\e_1|}\\
+&\[[\ol\pa,\pa]_G,\,\e_1]_G(-1)^{|\e_1|-1}(-1)^{|\e_1|}
\end{align}
Since $[\ol\pa, \pa]_G=\ol\pa\pa+\pa\ol\pa=0$, we have 
\bgn{align}
[\ol\pa, D\e_1]_G=&\[[\e_1,\ol\pa]_G,\,\pa]_G(-1)\\
=&\[[\ol\pa,\e_1]_G,\,\pa\]_G(-1)^{|\e_1|}\\
=&\[\pa,\,[\ol\pa, \e_1]_G\]_G(-1)^{|\e_1|}(-1)(-1)^{|\e_1|+1}\\
=&Dd_L\e_1
\end{align}
Substituting them, we obtain 
\bgn{align}
(-1)^{|\e_2|}d_L[\e_1, \e_2]_S
=&\[[\ol\pa, \e_2]_G,\,D\e_1\]_G(-1)^{|\e_2|\,(|\e_1|-1)}(-1)(-1)^{|\e_2|}
\\
+&\[Dd_L\e_1,\,\e_2\]_G(-1)^{|\e_2|}\\
\end{align}
Hence we have 
\bgn{align}
d_L[\e_1, \e_2]_S
=&\[D\e_1,\,[\ol\pa,\e_2]_G\]_G(-1)^{|\e_2|\,(|\e_1|-1)}
(-1)(-1)^{(|\e_1|-1)\,(|\e_2|+1)}\\
+&\[Dd_L\e_1,\,\e_2\]_G\\
=&\[Dd_L\e_1,\,\e_2\]_G+(-1)^{|\e_1|}\[D\e_1,\,[\ol\pa,\e_2]_G\]_G\\
=&[d_L\e_1,\e_2]_S+(-1)^{|\e_1|}[\e_1, d_L\e_2]_S
\end{align}
\qed\par
{\it Proof of lemma \ref{lem: S-3}}
For $\e_1, \e_2, \e_2\in \w^\bullet \ol{L}_\J$, it follow from (\ref{eq: app D}) that one have
\bgn{align*}
\big[\,\,[\e_1,\e_2]_S, \,\e_3\,\big]_S=&\big[\, D[ D\e_1, \e_2]_G, \, \e_3\,\big]_G\\
=&[\, [ D\e_1, D\e_2]_G, \, \e_3\,\,]_G(-1)^{(|\e_1|+1)} \\ \\
\big[\,\,[\e_2,\e_3]_S, \,\e_1\,\big]_S=&\big[\, D[ D\e_2, \e_3]_G, \, \e_1\,\big]_G\\
=&\big[\, [ D\e_2, \e_3]_G, \, D\e_1\,\big]_G (-1)^{(|\e_2|+|\e_3|)}\\\\
\big[\,\,[\e_3,\e_1]_S, \,\e_2\,\big]_S=&\big[\, D[ D\e_3, \e_1]_G, \, \e_2\,\big]_G\\
=&\big[\, D[ \e_3, D\e_1]_G, \, \e_2\,\big]_G(-1)^{(|\e_3|+1)}\\
=&\big[\, [ \e_3, D\e_1]_G, \, D\e_2\,\big]_G(-1)^{(|\e_3|+1)}(-1)^{(|\e_3|+|\e_1|)}
\end{align*}
Then we have three equations,
\bgn{align}
\big[\,\,[\e_1,\e_2]_S, \,\e_3\,\big]_S(-1)^{|\e_1||\e_3|}=&[\,\, [ D\e_1, D\e_2]_G, \, \e_3\,\,]_G(-1)^{(|\e_1|+1)} (-1)^{|\e_1||\e_3|} \\
\big[\,\,[\e_2,\e_3]_S, \,\e_1\,\big]_S(-1)^{|\e_2||\e_1|}=&\big[\, [ D\e_2, \e_3]_G, \, D\e_1\,\big]_G (-1)^{(|\e_2|+|\e_3|)}(-1)^{|\e_2||\e_1|}\\
\big[\,\,[\e_3,\e_1]_S, \,\e_2\,\big]_S(-1)^{|\e_3||\e_2|}=&\big[\, [ \e_3, D\e_1]_G, \, D\e_2\,\big]_G(-1)^{(|\e_3|+1)}(-1)^{(|\e_3|+|\e_1|)}(-1)^{|\e_3||\e_2|}
\end{align}
Multiplying $(-1)^{(|\e_3|-|\e_1|-1)}$, we have
\bgn{align}
(-1)^{(|\e_3|-|\e_1|-1)}\big[\,\,[\e_1,\e_2]_S, \,\e_3\,\big]_S(-1)^{|\e_1||\e_3|}=&[\,\, [ D\e_1, D\e_2]_G, \, \e_3\,\,]_G(-1)^{(|\e_1|+1)|\e_3|} \\
(-1)^{(|\e_3|-|\e_1|-1)}\big[\,\,[\e_2,\e_3]_S, \,\e_1\,\big]_S(-1)^{|\e_2||\e_1|}=&\big[\, [ D\e_2, \e_3]_G, \, D\e_1\,\big]_G(-1)^{(|\e_1|+1)(|\e_2|+1)} \\
(-1)^{(|\e_3|-|\e_1|-1)}\big[\,\,[\e_3,\e_1]_S, \,\e_2\,\big]_S(-1)^{|\e_3||\e_2|}=&\big[\, [ \e_3, D\e_1]_G, \, D\e_2\,\big]_G(-1)^{|\e_3|(|\e_2|+1)}
\end{align}

We apply the Jacobi identity of the graded bracket $[\, ,\, ]_G$ 
\bgn{align}
\big[\, [A,\,B]_G,\, C\big]_G (-1)^{|A||C|}
+\big[\,[ B, \, C]_G,\,A\big]_G(-1)^{|B||C|} +
\big[\,[C,\, A]_G,\, B\big]_G(-1)^{|C||B|}=0
\end{align}
Then we have the Jacobi identity of the Schouten bracket
$$
\big[\,\,[\e_1,\e_2]_S, \,\e_3\,\big]_S(-1)^{|\e_1||\e_3|}+
\big[\,\,[\e_2,\e_3]_S, \,\e_1\,\big]_S(-1)^{|\e_2||\e_1|}+
\big[\,\,[\e_3,\e_1]_S, \,\e_2\,\big]_S(-1)^{|\e_3||\e_2|}=0.
$$
\qed


\bgn{thebibliography}{99}
\bibitem{AN}
V.~Alexeev and V.~ Nikulin,
{\it Del Pezzo surfaces and K3 surfaces}, 
MSJ Memoirs, Mathematical Society of Japan, Vol. 15 (2006)

\bibitem{A.G.G}
V.~Apostolov, P.~Gauduchon, G.~Grantcharov,
{\it Bihermitian structures on complex surfaces},
Pro. London Math. Soc. {\bf 79}(1999), 414-428,
Corrigendum: {\bf 92}(2006), 200-202

\bibitem{Ch}
C.C.~Chevalley,
{\it The algebraic theory of Spinors}
Columbia University Press, 1954

\bibitem{De}
M.~Demazure, {\it Surfaces de del Pezzo II, III, IV, V}, 
 Lecture Notes in Math. 777, pp. 23-69, S\'eminaire sur les Singularit\'es des Surfaces, Palaiseau, France, 1976-1977
\bibitem{DV}
P.~Du Val,
{\it On isolated singularities of surfaces which do not affect the conditions of adjunction. I, II, III}, 
Proc. Cambridge Phi. Soc. 30 (1934), 453-465, 483-491

\bibitem{F}
A.~Fujiki and M.~ Pontecorvo,
Bihermitain anti-self-dual structures on Inoue surfaces, 
arXiv:0903.1320, to appear J.D.G
\bibitem{G.H.R}
S.~J. Gates C.~M. Hull and M. Ro\v cek,
{\it Twisted multiplets and new supersymmetric nonlinear $\sig$ models},
Nuclear Phys. B 248 (1984), 154-186
\bibitem{Go-1}
R.~Goto,
{\it Moduli spaces of topological calibrations,
Calabi-Yau, hyperK\"ahler, G$_2$ and Spin$(7)$ structures},
International Journal of Math. {\bf 115}, No. 3(2004), 211-257
\bibitem{Go0}
R.~Goto,
{\it On deformations of generalized Calabi-Yau, hyperK\"ahler, G$_2$ and Spin$(7)$ structures},
Math.DG/0512211
\bibitem{Go1}
R.~Goto,
{\it Deformations of \complex  and \K\"ahler structures},
Math. DG/0705.2495
\bibitem{Go2}
R.~Goto {\it Poisson structures and \K\"ahler submanifolds},
J. Math. Soc. Japan, vol. 61, No. 1 (2009), pp.107-132 
\bibitem{Gu1}
M.~Gualtieri,
{\it Generalized complex geometry},
Math.DG/0703298
\bibitem{Gu3}
M.~Gualtieri,
{\it Branes and Poisson varieties}
Math.DG/0710.2719
\bibitem{Hi3}
N.~Hitchin,
{\it Bihermitian metrics on Del Pezzo surfaces},
Math.DG/060821
\bibitem{Hu}
D.~Huybrechts,
{\it Generalized Calabi-Yau structures, K3 surfaces and B-fields}, math.AG/0306132,
International Journal of Math.16(2005) 13
\bibitem{Ko}
K.~Kodaira
{\it Complex manifolds and deformations of complex structures}
Grundlehren der Mathematischen Wissenschaften, {\bf 283}, springer-Verlag,
(1986)
\bibitem{KS}
K.~Kodaira and D.C.~Spencer,
{\it  On deformations of complex, analytic structures I,II}
Ann. of Math., 67(1958), 328-466
\bibitem{KSIII}
K.~Kodaira and D.C.~Spencer,
{\it On deformations of complex analytic structure, III.
stability theorems for complex structures}
Ann. of Math., 71(1960),43-76
\bibitem{KoSc}
Y.~Kosmann-Schwartzbach, 
{\it Derived brackets}, 
Lett. Math. Phys., 69, 61-87 (2004), 
math.DG/0312524

\bibitem{L.W.X}
Z.-J.~Liu, A.~Weinstein and Ping. ~Xu, 
{\it Manin triples for Lie bialgebroids},
J.Diff. Geom, {\bf 45},(1997) 547-574
\bibitem{Na}
Y.Namikawa
{\it Poisson deformations of affine symplectic varieties}
Math.AG/0609741

\bibitem{Pol}
A.~Polishchuk
{Algebraic geometry of Poisson varieties}
J. Math.Sci, Vol.84, No. 5, 1997, 1413-1444
\bibitem{Sa}
F.~Sakai,
{\it Anti-Kodaira dimension of ruled surfaces},
Sci. Rep. Saitama Univ. {\bf 2}(1982) 1-7
\bibitem{Sa1}
F.~Sakai, 
{\it Anticanonical Models of Rational Surfaces},
Math. Annalen, 269, 384-411 (1984)
\end{thebibliography}
\address{
Department of Mathematics \\ 
Graduate School of Science \\
Osaka University
Toyonaka, Osaka, 560\\
Japan
}
{goto@math.sci.osaka-u.ac.jp}
\end{document}